\PassOptionsToPackage{unicode}{hyperref}
\PassOptionsToPackage{hyphens}{url}
\documentclass[
  11pt]{article}
\usepackage{amsmath,amssymb}
\usepackage{iftex}
\ifPDFTeX
  \usepackage[T1]{fontenc}
  \usepackage[utf8]{inputenc}
  \usepackage{textcomp} 
\else 
  \usepackage{unicode-math} 
  \defaultfontfeatures{Scale=MatchLowercase}
  \defaultfontfeatures[\rmfamily]{Ligatures=TeX,Scale=1}
\fi
\usepackage{lmodern}
\ifPDFTeX\else
\fi
\IfFileExists{upquote.sty}{\usepackage{upquote}}{}
\IfFileExists{microtype.sty}{
  \usepackage[]{microtype}
  \UseMicrotypeSet[protrusion]{basicmath} 
}{}
\makeatletter
\@ifundefined{KOMAClassName}{
  \IfFileExists{parskip.sty}{%
    \usepackage{parskip}
  }{
    \setlength{\parindent}{0pt}
    \setlength{\parskip}{6pt plus 2pt minus 1pt}}
}{
  \KOMAoptions{parskip=half}}
\makeatother
\usepackage{xcolor}
\usepackage[margin=25mm]{geometry}
\usepackage{longtable,booktabs,array}
\usepackage{calc} 
\usepackage{etoolbox}
\makeatletter
\patchcmd\longtable{\par}{\if@noskipsec\mbox{}\fi\par}{}{}
\makeatother
\IfFileExists{footnotehyper.sty}{\usepackage{footnotehyper}}{\usepackage{footnote}}
\makesavenoteenv{longtable}
\usepackage{graphicx}
\makeatletter
\def\maxwidth{\ifdim\Gin@nat@width>\linewidth\linewidth\else\Gin@nat@width\fi}
\def\maxheight{\ifdim\Gin@nat@height>\textheight\textheight\else\Gin@nat@height\fi}
\makeatother
\setkeys{Gin}{width=\maxwidth,height=\maxheight,keepaspectratio}
\makeatletter
\def\fps@figure{htbp}
\makeatother
\usepackage{amsmath,amssymb,needspace}
\providecommand{\lt}{<}
\providecommand{\gt}{>}
\ifLuaTeX
  \usepackage{selnolig}  
\fi
\IfFileExists{bookmark.sty}{\usepackage{bookmark}}{\usepackage{hyperref}}
\IfFileExists{xurl.sty}{\usepackage{xurl}}{} 
\hypersetup{
  pdftitle={Planar Disk Filters and Real Chromatic Roots at Beraha Numbers},
  pdfauthor={Reza Mactabi},
  hidelinks,
  pdfcreator={LaTeX via pandoc}}

\title{Planar Disk Filters and Real Chromatic Roots at Beraha Numbers}
\author{Reza Mactabi\thanks{ORCID: \href{https://orcid.org/0009-0000-0925-2210}{0009-0000-0925-2210}. Independent researcher.}}
\date{3 October 2026}

\graphicspath{{output/beraha_submission/}{./}}
\begin{document}
\maketitle

\begin{abstract}
We construct simple plane triangulations whose
chromatic polynomials have real roots accumulating at every noninteger
Beraha number \(B_\ell=4\cos^2(\pi/\ell)\); the integer values
\(0,1,2,3\) are roots of \(K_4\). The construction rests on the
observation that a triangulated disk with four boundary vertices acts,
inside a larger triangulation, as a transfer operator lying in a
commutative three-parameter algebra. Rotating such disks by a quarter
turn and composing them in series moves their eigenvalue triples around
this algebra in a controllable way, and for every fixed
\(q\in(B_{10},4)\) one can manufacture a finite disk whose two competing
local eigenvalues are strictly dominated by the third --- a
\emph{physical filter}. Closing a grid of such disks between two poles
produces a sphere triangulation, and a fused Temperley--Lieb trace
identity accounts for its chromatic polynomial completely, sector by
sector, with no appeal to semisimplicity at the root of unity itself ---
generic semisimplicity enters only through polynomial identities
continued from \(\delta\gt2\) --- and no singular higher Jones--Wenzl
projectors. In the one-cap sector a Krein-type argument shows that every
sweep word has exactly one expanding eigenvalue; Gram-form inertia
controls every other sector uniformly in the width. At odd Beraha
indices the trace coefficient of the expanding eigenvalue vanishes; at
even indices two copies of that eigenvalue cancel between adjacent
sectors, and a short real-variable lemma converts either situation into
genuine real roots converging to the target. Indices below eleven are
treated by elementary families, exact outward-rounded integer
certificates, and one explicit disk at \(B_{10}\) combined with
Tutte\textquotesingle s golden identity. All graphs are ordinary,
simple, unweighted plane triangulations.
\end{abstract}
{\small MSC 2020: 05C31, 05C15, 82B20 · Keywords: chromatic polynomial, Beraha
numbers, Temperley--Lieb algebra, plane triangulations, transfer
matrices

}
\clearpage
{\small
\setcounter{tocdepth}{2}
\tableofcontents
}
\clearpage
\hypertarget{s1}{%
\section{The problem and the strategy}\label{s1}}

Write \(\chi_G(q)\) for the chromatic polynomial of a finite graph
\(G\). A \emph{simple plane triangulation} is a simple connected graph
embedded in the sphere with every face a triangle; choosing any face as
the outer face gives the equivalent plane formulation. The numbers

\[ B_\ell \;=\; 4\cos^2\!\frac{\pi}{\ell}, \qquad \ell = 2,3,4,\dots \]

are the Beraha numbers {[}Be{]}. We address the planar real-root
conjecture in the precise formulation of {[}JS, Conjecture 1.3{]}: for
every \(\ell\ge2\) and every \(\varepsilon\gt0\), a simple plane
triangulation has a real chromatic root within \(\varepsilon\) of
\(B_\ell\). The integer values are already roots of \(K_4\); the
substantive accumulation claim concerns the noninteger values.

{Theorem 1.1.} \emph{Every noninteger Beraha number --- \(B_\ell\) with
\(\ell=5\) or \(\ell\ge7\) --- is an accumulation point of real
chromatic roots of simple plane triangulations: for every
\(\varepsilon\gt0\) some simple plane triangulation has a real chromatic
root in
\((B_\ell-\varepsilon,\,B_\ell+\varepsilon)\setminus\{B_\ell\}\). At the
integer values \(B_2,B_3,B_4,B_6=0,1,2,3\), which are roots of
\(\chi_{K_4}\), no accumulation statement is made or needed.}

Two boundary remarks delimit what is and is not new in the statement.
Royle {[}R{]} proved that real chromatic roots of plane triangulations
approach \(4\); the endpoint is therefore not claimed here. Perrett and
Thomassen {[}PT{]} proved density of real chromatic roots for the much
wider class of planar graphs, away from a small interval containing
\(B_{10}\); their graphs are not triangulations, and the localization
demanded by Theorem~1.1 --- a root near each prescribed \(B_\ell\),
within the triangulation class --- does not follow from their result.
The quantum-group approach to Beraha accumulation predates this paper by
decades. Saleur {[}S{]} analyzed the relevant lattice systems and
explicitly stopped short of a full rigorous proof; Kauffman and Saleur
{[}KS{]} later described that accumulation argument as heuristic.
Jacobsen and Saleur {[}BCP{]} distinguished vanishing amplitudes from
sector crossings and explained exact finite-width coincidences through
representation theory --- mechanisms that are genuine antecedents of
what happens below. Biggs {[}B1,~B2{]} developed compatibility matrices
and representation-theoretic decompositions for cyclic graph families;
neither transfer matrices nor representation theory is claimed as new
here. Agol and Krushkal {[}AK{]} used Temperley--Lieb positivity to
bound chromatic \emph{values} of planar triangulations at Beraha
numbers; a bound on values does not by itself produce the real
\emph{roots} required by Theorem~1.1. The author\textquotesingle s
earlier manuscript {[}M{]} addresses particular finite targets and is
not a prerequisite here. The proposed contribution of this paper is the
disk construction at a variable target, the uniform all-width sector
analysis of the complete closure, and their assembly through real
sign-change lemmas into the full index range.

Theorem 1.1 concerns prescribed Beraha limit points. No claim is made
that these exhaust the possible real accumulation points of plane
triangulations.

\hypertarget{the-strategy-in-one-page}{%
\subsection{The strategy in one page}\label{the-strategy-in-one-page}}

It helps to see the whole arc before the details begin. The engine of
every transfer-matrix approach to this problem is the same: build a
one-parameter family of graphs \(G_t\) by repeating a gadget \(t\) times
around a cycle, so that \(\chi_{G_t}(q)\) becomes a finite sum
\(\sum_j \alpha_j(q)\,\lambda_j(q)^t\). If one eigenvalue \(\lambda_*\)
dominates all others near the target and its amplitude \(\alpha_*\)
changes sign simply at the target, then for large \(t\) the polynomial
must change sign nearby, and a real root appears and converges. Each
clause of that sentence hides a difficulty, and the proof must address
each of these difficulties.

First difficulty: \emph{dominance}. For raw gadgets the "entropic"
eigenvalues --- the ones counting bulk colorings --- tend to beat the
structured eigenvalue whose amplitude vanishes. Section 2 deals with
this head-on. We show that triangulated disks with four boundary
terminals form a commutative three-parameter calculus under serial
composition, with a quarter-turn rotation acting as an algebraic
involution-like move (Lemma 2.1 through display (D6)), and that inside
this calculus one can steer: for every fixed \(q\in(B_{10},4)\) there is
a finite disk whose eigenvalue triple \((u,v,e)\) satisfies
\(u\lt 0\lt v,e\lt -u\) (Proposition 2.2). Once each repeating cell is
such a \emph{filter}, both competing local channels contract and the
global dominance question can be addressed using the carrier and
Gram-form arguments of Sections 5 and 6; local contraction alone does
not settle it.

Second difficulty: \emph{completeness of the spectral accounting}. It is
fatally easy to analyze some sectors of a transfer operator and silently
postulate that no others contribute. Section 4 derives, from the
Fendley--Krushkal chromatic algebra map and the standard quantum trace,
a single boxed identity (T6) expressing \(\chi_{G_t}\) as an explicitly
weighted sum of traces over \emph{fused} link-state sectors, with bases
constructed directly and no division by vanishing quantum integers. The
identity is checked symbolically on nine complete graph polynomials in
Appendix B, and its derivation is diagrammatic.

Third difficulty: \emph{what happens at the root of unity}. At
\(q=B_\ell\) the Temperley--Lieb algebra is not semisimple, and
standard-module Gram forms degenerate. Rather than fighting this, the
proof uses it. Section 5 shows that the one-cap sector carries, for
every sweep word, exactly one expanding eigenvalue --- a Krein-space
argument on a tridiagonal form \(H\) (Lemma 5.1), with the eigenvalue
given in closed form (H4) for the ascending sweep. Section 6 then
controls every other sector by Gram inertia: positive-definite sectors
give operator norms at most one; the single semidefinite sector has a
one-dimensional radical on which everything acts trivially; and at even
indices the celebrated reflection structure of Ridout and Saint-Aubin
produces a \emph{second copy} of the carrier inside the adjacent sector,
whose quotient we prove positive by a derivative argument (Proposition
6.2, displays (G6a)--(G6b)).

Fourth difficulty: \emph{real roots rather than complex ones}. The
Beraha--Kahane--Weiss theory {[}BKW, BO{]} produces limit points of
complex zeros; it cannot by itself produce real zeros, and the
distinction is not pedantic --- \(f_m(x)=x^2+1/m\) has complex zeros
converging to the real point \(0\) and no real zero at all. Section 7
therefore proves two self-contained sign-change lemmas tailored to the
exact hypotheses established earlier: one for a single vanishing
amplitude (the odd-index mechanism), one for a cancelling pair of
amplitudes attached to branches that may coincide at the target (the
even-index mechanism, where no slope-separation assumption is available
or needed).

The assembly is then short. Section 8 proves Theorem 1.1 for every
\(\ell\ge 11\), both parities. Section 9 handles \(2\le\ell\le 10\): the
integers via \(K_4\), \(B_5\) via bipyramids, \(B_7,B_8,B_9\) via the
elementary disk and exact interval certificates, and the golden boundary
\(B_{10}\) via one explicit composite disk with \(419{,}495\) internal
vertices. Section 10 proves the produced roots are distinct from their
targets --- by a conjugate-in-\((0,1)\) argument at most indices, and by
Tutte\textquotesingle s golden identity at \(B_{10}\). Section 11
explains precisely which statements rest on the executable certificates
and which are proved for all widths at once.

\Needspace{16\baselineskip}A notational table, used throughout:

\begin{longtable}[]{@{}p{\dimexpr.42\linewidth-2\tabcolsep}p{\dimexpr.58\linewidth-2\tabcolsep}@{}}
\toprule\noalign{}
Symbol & Meaning \\
\midrule\noalign{}
\endhead
\bottomrule\noalign{}
\endlastfoot
\(\ell\) & Beraha index; the target is \(q_*=B_\ell\) \\
\(n\) & number of internal sites of a path between two fixed poles \\
\(N=2n+2\) & number of Temperley--Lieb strands \\
\(t\) & number of sweeps in the closed graph; we always take
\(t=2m\ge4\) \\
\(a=q-2,\ \tau=q-3\) & local color parameters \\
\((u,v,e)\) & the three local eigenvalues of a disk \\
\(g=-u,\ r_v=v/g,\ r_e=e/g\) & normalization used when \(u\lt 0\) \\
\(W_j=W_{N,2j},\ V_j=FW_j\) & standard and fused link-state modules \\
\(\mathsf B,\mathsf C,D_0\) & disk parameters (sans serif, so no clash
with the Beraha numbers \(B_\ell\)) \\
\(A(q),B(q)\) & trace amplitudes, always written with their argument;
Beraha numbers always carry a subscript \\
\(\Gamma,\Xi,\Theta\) & blocks of the degenerating Gram matrix in
(G6) \\
\end{longtable}

\hypertarget{s2}{%
\section{Four-terminal disks and physical filters}\label{s2}}

\hypertarget{the-local-three-dimensional-algebra}{%
\subsection{The local three-dimensional
algebra}\label{the-local-three-dimensional-algebra}}

Let \(D\) be a graph embedded in a disk with boundary cycle
\((l,x,r,X)\), where boundary-adjacent colors are required to differ. We
think of \(l\) and \(r\) as \emph{poles} and of \(x,X\) as the input and
output of a path that the disk will join. For a fixed proper boundary
coloring, let \(W_D\) count its extensions to proper colorings of the
internal vertices. The first lemma says that this count, as a function
of the boundary coloring, lives in a three-dimensional space --- and the
reason is planarity, not any symmetry we impose.

{Lemma 2.1.} \emph{There are polynomials \(u,\mathsf B,D_0\) in \(q\)
such that, on proper boundary assignments,}
\[ W_D \;=\; D_0 \;+\; u\,\mathbf 1_{x=X} \;+\; (\mathsf B-D_0)\,\mathbf 1_{l=r}. \tag{D1} \]

{Proof.} Expand all nonboundary edge constraints, including any boundary
chords, by deletion--contraction with the boundary colors held fixed.
Each resulting term induces a noncrossing connectivity partition of the
four boundary vertices, weighted by a polynomial coming from the
components that contain no boundary vertex. A partition identifying two
boundary-adjacent vertices contributes nothing on proper boundary
assignments. Two separate blocks joining \(l\) to \(r\) and \(x\) to
\(X\) cannot both occur, because inside a disk those two connections
would have to cross; and a single block containing all four terminals
again identifies adjacent colors. The only surviving equality conditions
are therefore the empty one and the two opposite-pair identifications
taken separately. This proves (D1) for all sufficiently large integer
\(q\), hence as a polynomial identity. {\(\square\)}

To convert counts into operators, let \(U_i\) denote the map that sums
over the color at site \(i\) of a path, holding all other colors fixed
and preserving properness. On a fiber where the two neighbors of site
\(i\) are colored unequally, \(U_i\) is an all-ones matrix of size
\(a=q-2\); when the neighbors agree it has size \(a+1\). Consequently

\[ U_i^2-aU_i=\mathbf 1_{l=r}\,U_i, \qquad U_i\,(U_i-aI)\,(U_i-(a+1)I)=0, \tag{D2} \]

and the transfer operator of the disk \(D\), acting at site \(i\), is
the specific polynomial

\[ T_i \;=\; uI \;+\; D_0\,U_i \;+\; (\mathsf B-D_0)\,(U_i^2-aU_i). \tag{D3} \]

Its three local eigenvalues --- the data we will steer --- are

\[ u,\qquad v=u+aD_0,\qquad e=u+(a+1)\mathsf B. \tag{D4} \]

A fact worth pausing on: the two zero-sum fibers share the same
eigenvalue \(u\). This is a consequence of planarity, visible already in
(D1), and not an extra symmetry assumption. Writing \(\mathsf C=D_0+u\),
the parameters invert to

\[ \mathsf B=\frac{e-u}{a+1},\qquad D_0=\frac{v-u}{a},\qquad \mathsf C=\frac{v+(a-1)u}{a}, \tag{D5} \]

and the denominators are harmless at our noninteger targets, all of
which exceed \(2\). For physical disks the polynomial meaning of these
quotients also follows directly from (D1).

\hypertarget{rotation-and-serial-composition}{%
\subsection{Rotation and serial
composition}\label{rotation-and-serial-composition}}

Two operations generate everything we need. A \emph{quarter turn} of the
entire disk interchanges the roles of the two opposite boundary pairs;
on parameters this swaps \(\mathsf B\) and \(\mathsf C\) while fixing
\(D_0\), and on eigenvalue triples it acts by

\[ \mathcal R(u,v,e) \;=\; \bigl(\mathsf B-D_0,\;\; \mathsf B+(a-1)D_0,\;\; \mathsf B-D_0+(a+1)\mathsf C\bigr). \tag{D6} \]

\emph{Serial composition} glues the output path \((l,X,r)\) of one disk
to the input path of the next; summing over the shared middle color
multiplies the transfer operators, and since both lie in the commutative
algebra generated by (D3), their eigenvalue triples simply multiply
componentwise. The calculus is thus a semigroup of triples, acted on by
the rotation \(\mathcal R\).

The basic generators are the \emph{fans}. Let \(F_k\) consist of a path
of \(k\) edges from \(x\) to \(X\), with both poles adjacent to every
path vertex. It is a triangulated disk with \(k-1\) internal vertices,
and a one-line induction from \(F_1\) (whose transfer is \(U-I\)) gives
the clean triple

\[ F_k:\qquad \bigl((-1)^k,\;\tau^k,\;a^k\bigr). \tag{D7} \]

\begin{figure}
\centering
\includegraphics[width=6.25in,height=\textheight]{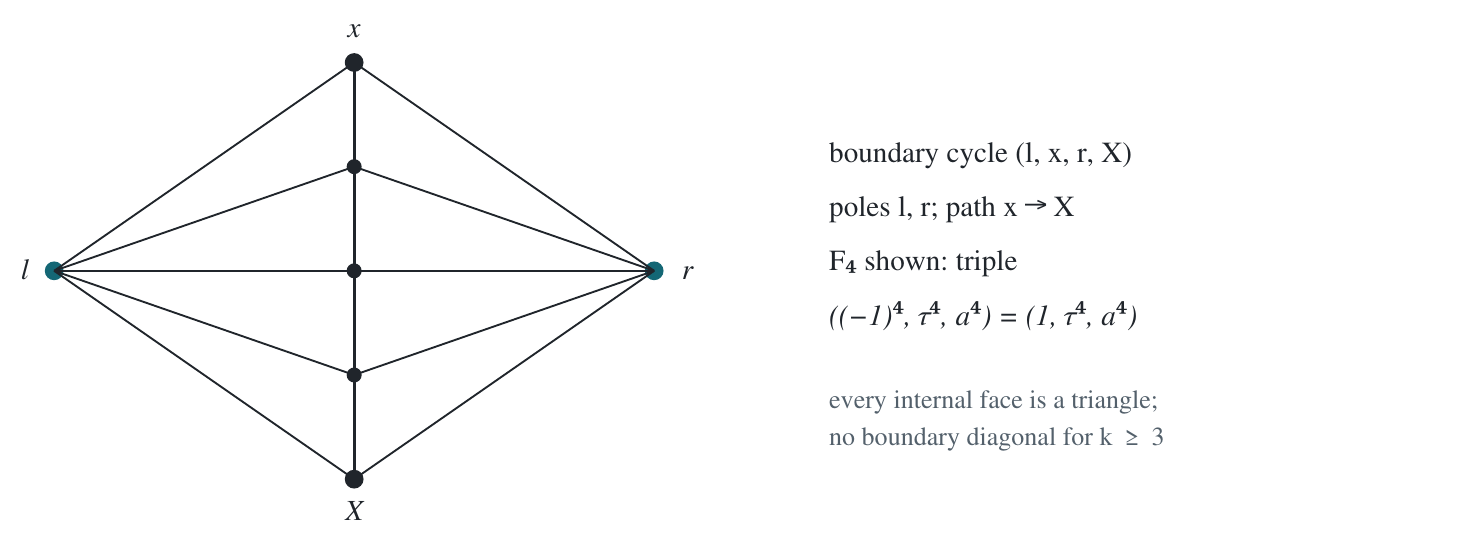}
\caption{The fan \(F_k\) (here \(k=4\)): a path of \(k\) edges from
\(x\) to \(X\) with both poles adjacent to every path vertex. Its
eigenvalue triple is \(((-1)^k,\tau^k,a^k)\) by (D7); the boundary cycle
\((l,x,r,X)\) is the interface along which disks are rotated, composed,
and glued into faces.}
\end{figure}

For \(k\ge 3\) the fan has neither boundary diagonal; rotation preserves
this property; and serial composition of two simple disks without
boundary diagonals is again a simple triangulated disk without boundary
diagonals, since the overlap of the two is exactly the identified
two-edge boundary path and no edge joins its nonconsecutive vertices. If
the internal vertex counts are \(r_1\) and \(r_2\), the composite has
\(r_1+r_2+1\). These bookkeeping remarks matter later: they are what
guarantees that the elaborate composites we build remain legal cells for
a simple triangulation.

\begin{figure}
\centering
\includegraphics[width=6.25in,height=\textheight]{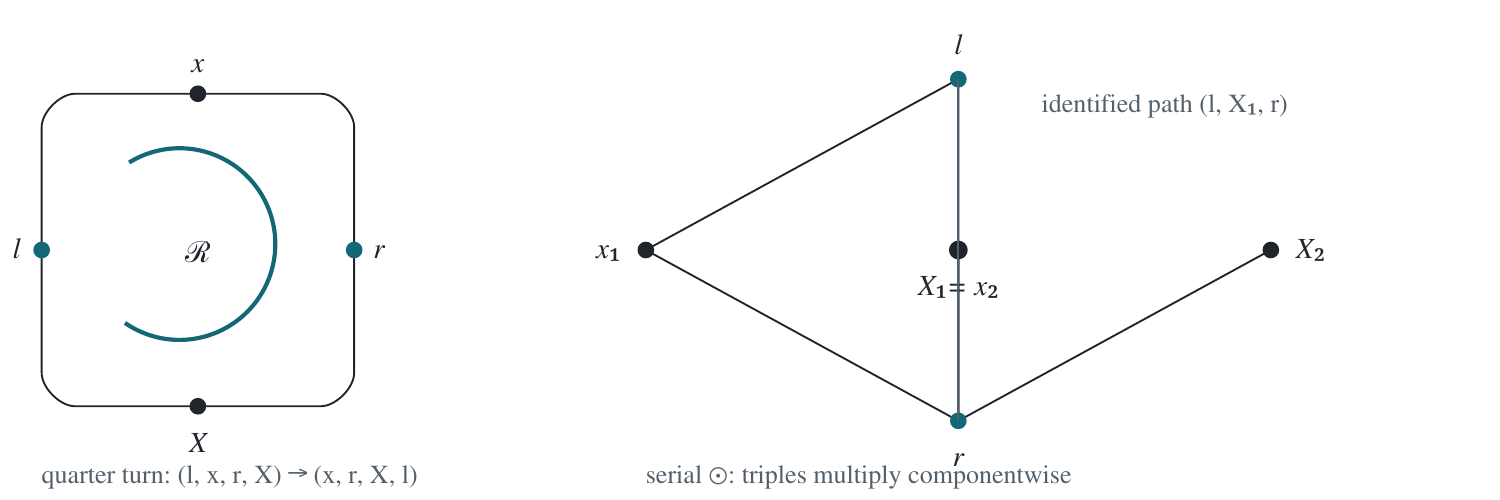}
\caption{The two moves of the disk calculus. Left: the quarter turn
\(\mathcal R\), which exchanges the roles of the opposite boundary pairs
--- on parameters, \(\mathsf B\leftrightarrow\mathsf C\) with \(D_0\)
fixed, display (D6). Right: serial composition \(\odot\) along the path
axis: the output path \((l,X_1,r)\) of the first disk is identified with
the input path \((l,x_2,r)\) of the second (dashed seam), the poles
\(l,r\) being shared; the composite has boundary \((l,x_1,r,X_2)\), and
transfer operators multiply in the commutative algebra, so eigenvalue
triples multiply componentwise.}
\end{figure}

\hypertarget{existence-of-a-strict-filter-above-b_10}{%
\subsection{\texorpdfstring{Existence of a strict filter above
\(B_{10}\)}{Existence of a strict filter above B\_\{10\}}}\label{existence-of-a-strict-filter-above-b_10}}

Call a disk a \emph{filter} at \(q\) if its triple satisfies

\[ u \lt 0, \qquad 0 \lt v,\,e \lt -u. \tag{D8} \]

Under (D8) the normalized operator \(-T/|u|\) has local spectrum
\(\{1,-r_v,-r_e\}\) with \(0\lt r_v,r_e\lt 1\): both competing channels
contract, strictly, before any global analysis begins.

{Proposition 2.2.} \emph{For each fixed \(q\in(B_{10},4)\) there exists
a finite simple triangulated disk without boundary diagonals satisfying
(D8). It may be taken of the form
\(\mathcal R\!\bigl(R_k^{\,\odot m_1}\bigr)\odot R_k^{\,\odot m_2}\)
with \(R_k=\mathcal R(F_k)\), \(k\) odd, where \(\odot\) denotes serial
composition.}

The proof is organized so that no two limits interact: a
closure-and-openness argument (Step 4) reduces everything to a single
explicit limit (Step 2) and a finite calculus verification at one point
(Step 3). Throughout, fix \(q\in(B_{10},4)\) and set

\[ a=q-2\in(\phi,2),\qquad \tau=a-1\in(\phi-1,1),\qquad \kappa=a-\frac1a\in(1,\tfrac32); \]

the hypothesis \(q\gt B_{10}=\phi+2\) is used once, and only through
\(\kappa\gt1\).

\Needspace{5\baselineskip}\textbf{Step 0: the projective semigroup of legal disks.} Call a disk
\emph{legal} if it is simple, triangulated, and free of boundary
diagonals. Since the triple of a disk determines its transfer only up to
the data in (D4), and (D8) is invariant under multiplication of the
triple by a positive scalar, it is natural to work in the space
\(\mathbb S=(\mathbb R^3\setminus\{0\})/\mathbb R_{\gt0}\) of positive
rays, with the quotient topology. Serial composition descends to a
continuous componentwise product \(\odot\) on the open subset of
\(\mathbb S\times\mathbb S\) where the product triple is nonzero (the
physical triple of a composite \emph{is} the componentwise product, by
the commutative-algebra computation of Section 2.2, so no scale is ever
discarded), and the rotation \(\mathcal R\) of (D6) is \emph{linear} in
\((u,v,e)\) --- the parameters \(\mathsf B,\mathsf C,D_0\) of (D5) are
linear functions of the triple --- so it descends to a continuous map
wherever its image is nonzero. Let

\[ \mathcal T\;=\;\{\,[\text{triple of }D]\in\mathbb S:\ D\ \text{a legal disk with nonzero triple}\,\},\qquad \overline{\mathcal T}=\text{its closure in }\mathbb S. \]

{Lemma 2.2a (legality and closure properties).} \emph{(i) For odd
\(k\ge3\), \(F_k\) is legal, and legality is preserved by rotation and
by serial composition. (ii) \(\mathcal T\) is closed under \(\odot\) and
under \(\mathcal R\) whenever the resulting triple is nonzero;
consequently so is \(\overline{\mathcal T}\), at every point where the
relevant image is nonzero.}

{Proof.} (i) \(F_k\) for \(k\ge3\) is a simple triangulated disk, and
its only candidate boundary diagonals \(lr\), \(xX\) are absent since
\(l,r\) are nonadjacent and the \(x\)--\(X\) path has length \(k\ge3\).
Rotation relabels the boundary and changes nothing internal. For a
serial composite of legal disks, the two copies intersect exactly in the
identified boundary path \((l,X_1{=}x_2,r)\); no edge joins
nonconsecutive vertices of that path in either factor (no boundary
diagonals), so the union is simple; the glued region is triangulated by
its parts; and a boundary diagonal of the composite would be an edge
\(lr\) (absent in both factors) or \(x_1X_2\) (impossible: \(x_1,X_2\)
lie in different factors and are not on the shared path). (ii) The first
claim restates Section 2.2 plus part (i); the second follows from
continuity of \(\odot\) and \(\mathcal R\): if \(x_m\to x\),
\(y_m\to y\) with \(x_m,y_m\in\mathcal T\), then
\(x_m\odot y_m\to x\odot y\), and \(\mathcal R(x_m)\to\mathcal R(x)\)
whenever \(\mathcal R(x)\ne0\). {\(\square\)}

\Needspace{5\baselineskip}\textbf{Step 1: the rotated fans.} For odd \(k\), the triple of
\(R_k=\mathcal R(F_k)\) is, from (D6) and (D7),

\[ u_k=\frac{a^{k+1}-(a+1)\tau^k-1}{a(a+1)},\qquad v_k=u_k+1+\tau^k,\qquad e_k=u_k+\frac{a+1}{a}\,\bigl(\tau^k-\tau\bigr). \tag{D9} \]

Since \(a\gt\phi\) and \(0\lt\tau\lt1\),

\[ u_k\;\ge\;\frac{a^{k+1}-a-2}{a(a+1)}\;\gt\;0\quad(k\ge3),\qquad u_k\;\xrightarrow[k\to\infty]{}\;\infty\ \ \text{monotonically in this bound}. \tag{D9a} \]

Define \(b_k=v_k/u_k\gt1\) and \(c_k=e_k/u_k\). Because
\(\tau^k\lt\tau\), we have \(c_k\lt1\); and \(c_k\gt0\) as soon as
\(u_k\gt\tfrac{(a+1)\tau}{a}\), which (D9a) guarantees for all odd
\(k\ge K_1(q)\), a threshold determined by (D9a). Finally, with
\(x_k=(1+\tau^k)/u_k\) and \(y_k=\tfrac{(a+1)(\tau-\tau^k)}{a\,u_k}\),
the inequality \(|\log(1+x)-x|\le x^2\) for \(|x|\le\tfrac12\) gives,
for all odd \(k\ge K_2(q)\),

\[ \begin{gathered}\Bigl|\log b_k-x_k\Bigr|\le x_k^2,\qquad \Bigl|\log c_k+y_k\Bigr|\le y_k^2,\\ \frac{y_k}{x_k}=\frac{(a+1)(\tau-\tau^k)}{a(1+\tau^k)}\;\xrightarrow[k\to\infty]{}\;\frac{(a+1)\tau}{a}\cdot 1=\frac{a^2-1}{a}=\kappa.\end{gathered} \tag{D9b} \]

so that \(\log b_k\to0^+\) and \(\log c_k/\log b_k\to-\kappa\), with all
error terms explicit in \(u_k\) and \(\tau^k\).

\Needspace{5\baselineskip}\textbf{Step 2: one limit --- the curve lies in the closure.}

{Lemma 2.2b.} \emph{For every \(s\gt0\), the ray of
\(X(s)=\bigl(1,\,e^{s},\,e^{-\kappa s}\bigr)\) belongs to
\(\overline{\mathcal T}\).}

{Proof.} For odd \(k\ge\max(K_1,K_2)\) set
\(m_k=\bigl\lceil s/\log b_k\bigr\rceil\) and consider the legal disk
\(R_k^{\odot m_k}\), whose ray is that of
\((1,\,b_k^{m_k},\,c_k^{m_k})\) after dividing by \(u_k^{m_k}\gt0\).
Then
\[ \bigl|m_k\log b_k-s\bigr|\le\log b_k\xrightarrow[k\to\infty]{}0, \qquad m_k\log c_k=\bigl(m_k\log b_k\bigr)\cdot\frac{\log c_k}{\log b_k}\xrightarrow[k\to\infty]{}s\cdot(-\kappa), \]
the first by the definition of the ceiling, the second by (D9b). Hence
\((1,b_k^{m_k},c_k^{m_k})\to X(s)\) in \(\mathbb R^3\), and rays
converge to rays. {\(\square\)}

\Needspace{5\baselineskip}\textbf{Step 3: a finite calculus verification at one point.}

{Lemma 2.2c.} \emph{There exist \(s^*\gt0\) and \(h\gt0\), depending
only on \((a,\kappa)\), such that the triple
\(Y=\mathcal R\bigl(X(s^*)\bigr)\odot X(h)\) satisfies the strict
inequalities (D8).}

{Proof.} Write \(\mathcal R(X(s))=(u_s,v_s,e_s)\). From (D5)--(D6) with
\((u,v,e)=(1,e^s,e^{-\kappa s})\):
\[ u_s=\frac{e^{-\kappa s}-1}{a+1}-\frac{e^{s}-1}{a},\qquad v_s=\frac{e^{-\kappa s}-1}{a+1}+(a-1)\,\frac{e^{s}-1}{a},\qquad e_s=u_s+(a+1)\Bigl(\frac{e^{s}-1}{a}+1\Bigr). \tag{D10} \]
These are exact. Expanding
\(e^{\pm\sigma s}=1\pm\sigma s+\tfrac{\sigma^2s^2}{2}+\theta\),
\(|\theta|\le\tfrac{\sigma^3s^3}{6}e^{\sigma s}\), the first-order
coefficient of \(u_s\) is
\[ -\Bigl(\frac{\kappa}{a+1}+\frac1a\Bigr)=-\frac{\kappa a+a+1}{a(a+1)}=-\frac{(a^2-1)+a+1}{a(a+1)}=-1, \]
and the first-order coefficient of \(v_s\) is
\(-\tfrac{\kappa}{a+1}+\tfrac{a-1}{a}=\tfrac{-(a^2-1)+(a^2-1)}{a(a+1)}=0\),
with second-order coefficient
\[ c_0=\frac{\kappa^2}{2(a+1)}+\frac{a-1}{2a}=\frac{\kappa^2a+a^2-1}{2a(a+1)}=\frac{\kappa(\kappa+1)}{2(a+1)}\;\gt\;0 . \]
Therefore, with all remainders bounded by the displayed Lagrange terms,
there are constants \(s_0(a)\gt0\) and \(C_1(a)\) such that for
\(0\lt s\le s_0\):
\[ -\tfrac32 s\le u_s\le-\tfrac12 s,\qquad \tfrac12 c_0s^2\le v_s\le\tfrac32 c_0s^2,\qquad a+1-C_1s\le e_s\le a+1+C_1s . \tag{D10a} \]
Shrink \(s_0\) so that \(C_1s_0\lt(a+1)/2\). Set \(V_s=v_s/(-u_s)\) and
\(E_s=e_s/(-u_s)\); then by (D10a)
\[ \tfrac13 c_0 s\;\le\;V_s\;\le\;3c_0 s,\qquad \tfrac{a+1-C_1s}{3s/2}\;\le\;E_s\;\le\;\tfrac{2(a+1+C_1s)}{s}, \]
whence
\(E_sV_s^{\kappa}\le 2(a+1+C_1s_0)\,(3c_0)^{\kappa}\,s^{\kappa-1}\).
Because \(\kappa\gt1\) --- the sole use of \(q\gt B_{10}\) --- the right
side tends to \(0\); choose \(s^*\in(0,s_0]\) so that
\(E_{s^*}V_{s^*}^{\kappa}\lt1\) and \(E_{s^*}\gt V_{s^*}\) (the latter
holds for all small \(s\) since \(E_s\to\infty\), \(V_s\to0\)), and set
\(h=\log(E_{s^*}/V_{s^*})/(\kappa+1)\gt0\). The componentwise product
\(Y=\mathcal R(X(s^*))\odot X(h)\) has
\[ u(Y)=u_{s^*}\lt0,\qquad \frac{v(Y)}{-u(Y)}=V_{s^*}e^{h}=\bigl(E_{s^*}V_{s^*}^{\kappa}\bigr)^{\frac1{\kappa+1}}=E_{s^*}e^{-\kappa h}=\frac{e(Y)}{-u(Y)}, \]
a common value in \((0,1)\), and both numerators positive: (D8) holds
strictly at \(Y\). {\(\square\)}

\Needspace{5\baselineskip}\textbf{Step 4: assembly.}

{Proof of Proposition 2.2.} By Lemma 2.2b, the rays of \(X(s^*)\) and
\(X(h)\) lie in \(\overline{\mathcal T}\);
\(u(\mathcal R(X(s^*)))=u_{s^*}\ne0\), so by Lemma 2.2a(ii) the ray of
\(Y=\mathcal R(X(s^*))\odot X(h)\) lies in \(\overline{\mathcal T}\).
The set of rays satisfying the three strict inequalities (D8) is open in
\(\mathbb S\) (in the chart \(u=-1\) it is the open set
\(\{0\lt v\lt1,\ 0\lt e\lt1\}\)), and it contains the ray of \(Y\) by
Lemma 2.2c. An open neighborhood of a closure point meets the set
itself: some legal disk has its ray in that neighborhood, and therefore
satisfies (D8). Unwinding Lemmas 2.2a--2.2b, that disk may be taken of
the stated form \(\mathcal R(R_k^{\odot m_1})\odot R_k^{\odot m_2}\) for
a single sufficiently large odd \(k\), with
\(m_1=\lceil s^*/\log b_k\rceil\), \(m_2=\lceil h/\log b_k\rceil\):
convergence of each factor and continuity of \(\mathcal R\) and
\(\odot\) make the composite rays converge to \(Y\), and openness
finishes as before. The thresholds \(K_1,K_2,s_0,C_1,s^*\) are effective
--- each is determined by the displayed inequalities, without asserting
that these thresholds were numerically evaluated in the supplied audit
--- but no closed-form expressions and no neighborhood sizes are
claimed, and none are needed: the argument is existential, and the
certificates of Section 9 are where explicitness is delivered.
{\(\square\)}

{Remark (fixed-disk discipline).} The disk produced depends on the
target \(q_*\) and is then held fixed: its triple \((u,v,e)\) consists
of polynomials in \(q\), the inequalities (D8) hold at \(q_*\), hence on
a neighborhood by continuity, and the scale \(g=-u\) entering (T6) is
the true polynomial eigenvalue of this fixed disk. The limiting
construction is never allowed to vary the graph with \(q\) inside a
sign-change argument. The explicit disk of Section 9.4 at the boundary
value \(B_{10}\), where \(\kappa=1\) and this proposition does not
apply, is built from the same two-stage shape (S4) with the parameters
chosen by exact arithmetic instead of the limit.

\hypertarget{s3}{%
\section{Closing the disks into simple triangulations}\label{s3}}

Fix \(n\ge2\) and \(t\ge3\). Take two poles \(H,O\) and vertices
\(v_{i,r}\) for \(1\le i\le n\) and \(r\in\mathbb Z/t\mathbb Z\), with
the conventions \(v_{0,r}=H\) and \(v_{n+1,r}=O\). Define oriented
quadrilateral faces

\[ Q_{i,r}=\bigl(v_{i-1,r+1},\,v_{i,r},\,v_{i+1,r},\,v_{i,r+1}\bigr). \tag{Q1} \]

{Lemma 3.1.} \emph{These faces form a simple quadrangulation of the
sphere with \(nt+2\) vertices, \(2nt\) edges, and \(nt\) faces.}

{Proof.} Embed the \(t\) paths \(H,v_{1,r},\dots,v_{n,r},O\) as
consecutive meridians with disjoint interiors. The region between
meridians \(r\) and \(r+1\) is a disk with \(2n+2\) boundary edges;
adding the noncrossing diagonals \(v_{i+1,r}v_{i,r+1}\) for
\(1\le i\lt n\) cuts it into exactly the faces (Q1). For \(t\ge3\) no
two added edges coincide, and none coincides with a path edge. The
counts follow. {\(\square\)}

\begin{figure}
\centering
\includegraphics[width=6.25in,height=\textheight]{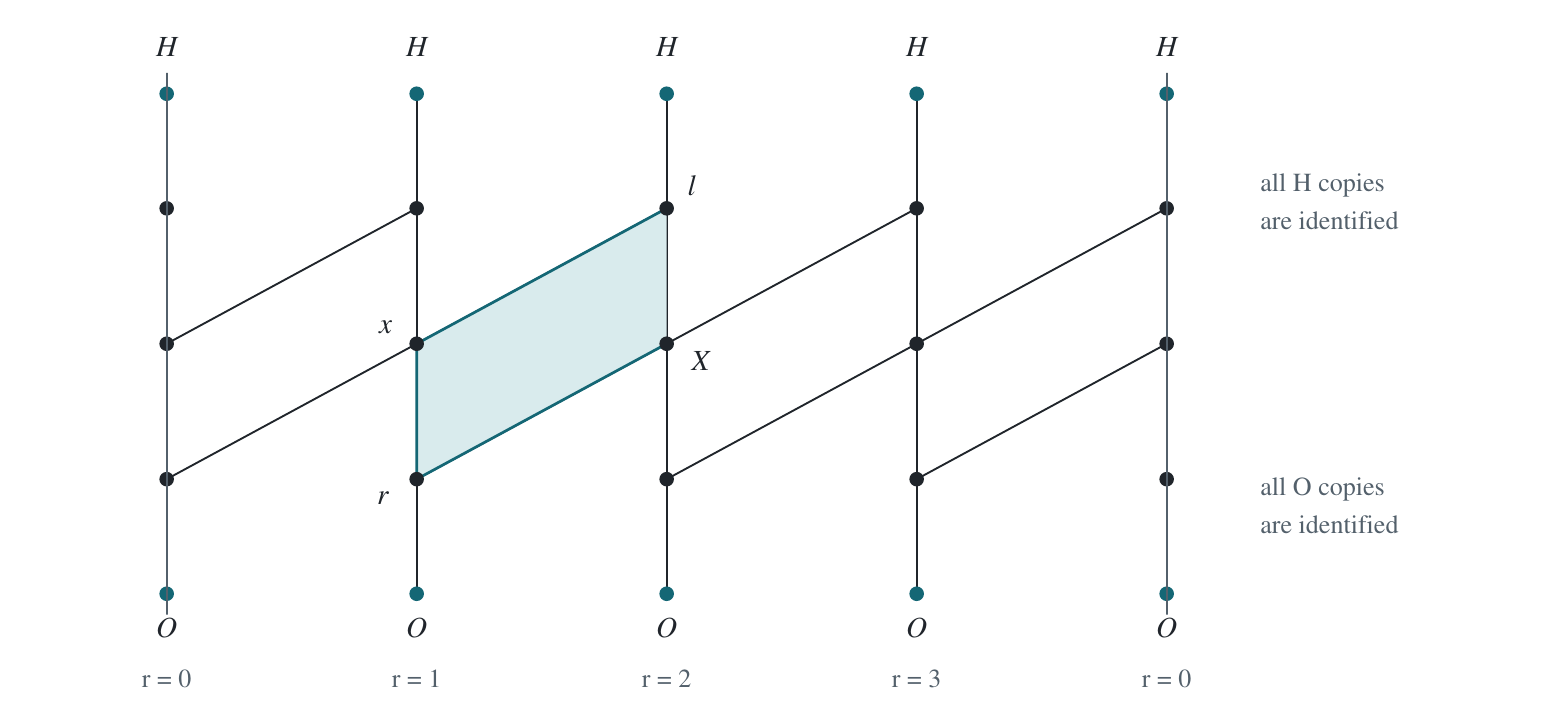}
\caption{An unrolled view of the closure with \(n=3,t=4\). Identify the
two dashed meridians, all copies of \(H\), and all copies of \(O\). The
shaded face is \(Q_{2,1}=(v_{1,2},v_{2,1},v_{3,1},v_{2,2})\), labeled
\((l,x,r,X)\). The duplicated drawing vertices at the seam and poles
denote identifications, not extra graph vertices.}
\end{figure}

Now replace every face by a copy of a fixed disk \(D\), matching
terminals in the displayed cyclic order, with the internal vertices of
different copies disjoint. If \(D\) is simple and has no boundary
diagonals, the result \(G_t(D,n)\) is a simple sphere triangulation; if
\(r(D)\) is the number of internal vertices of \(D\), then

\[ |V(G_t)| \;=\; 2+nt\,\bigl(r(D)+1\bigr). \tag{Q2} \]

The elementary disk \(F_1\) is also a permitted special case: its single
diagonal \(xX\) becomes \(v_{i,r}v_{i,r+1}\), and for \(t\ge3\) these
diagonals are distinct and absent from the quadrangulation, so the
closure is again simple.

Dynamically, a \emph{sweep} updates the sites \(1,\dots,n\) in order.
When site \(i\) is updated, its left neighbor has already been updated
while its right neighbor has not, so the local face is exactly (Q1) and
a sweep realizes the product \(T_n\cdots T_1\). Cyclic identification of
the path states closes the trace. We always take \(t=2m\ge4\), which
removes a global sign from the normalization in the next section.

\hypertarget{s4}{%
\section{The complete fused trace}\label{s4}}

This section proves the boxed trace identity on which the rest of the
paper runs. It has three independent ingredients, proved separately: a
purely combinatorial \emph{transfer identity} expressing \(\chi_{G_t}\)
as a trace of operators on colorings (Lemma 4.1); an explicit
\emph{representation} of the Temperley--Lieb algebra on the coloring
space, under which the disk transfer operators of Section 2 are
literally the images of Temperley--Lieb elements (Lemma 4.2); and a
\emph{trace evaluation} identifying the coloring trace of every diagram
with the standard quantum trace (Lemma 4.3). The boxed identity (T6)
then follows by composing the three with the generic sector
decomposition (T4). Nothing here depends on semisimplicity at roots of
unity or on singular higher Jones--Wenzl projectors; the
chromatic-algebra map of Fendley and Krushkal {[}FK1{]} is the
conceptual source of Lemma 4.2, and the elementary identities below are
proved explicitly; the generic standard-module structure and
Markov-trace characterization are classical inputs {[}RS, LK{]}. The
factors of \(F\) are the nonsingular rank-two Jones--Wenzl idempotents,
which divide only by \(\delta\).

\hypertarget{the-transfer-identity-on-colorings}{%
\subsection{The transfer identity on
colorings}\label{the-transfer-identity-on-colorings}}

Fix \(n\ge2\), \(t\ge3\), and for each face position \((i,r)\),
\(1\le i\le n\), \(r\in\mathbb Z/t\), a triangulated disk \(D_{i,r}\)
with boundary cycle \((l,x,r,X)\), no boundary diagonals unless the disk
is \(F_1\), and internal vertex sets pairwise disjoint. Let \(G_t\) be
the closure of Section 3 with \(D_{i,r}\) glued into the face
\(Q_{i,r}\) of (Q1) in the displayed cyclic order. Write
\(\mathcal P(c_H,c_O)\) for the set of proper colorings of the path
\(H,x_1,\dots,x_n,O\) --- properness meaning adjacent entries differ
along the \(n+1\) path edges --- with the pole colors frozen at
\((c_H,c_O)\), and let \(\mathbb C\,\mathcal P(c_H,c_O)\) be its free
vector space. For a disk \(D\) placed at site \(i\), define \(T_i^{D}\)
on this space by

\[ \bigl(T_i^{D}\,f\bigr)(\dots,x_i',\dots)\;=\;\sum_{x_i}\,W_{D}\bigl(l{=}x_{i-1},\,x{=}x_i,\,r{=}x_{i+1},\,X{=}x_i'\bigr)\,f(\dots,x_i,\dots), \]

with \(x_0=c_H\), \(x_{n+1}=c_O\), and \(W_D\) the boundary-to-count
function of Section 2, which vanishes unless its four arguments form a
proper coloring of the boundary cycle.

{Lemma 4.1 (transfer identity).} \emph{With
\(T^{(r)}=T_n^{D_{n,r}}\cdots T_1^{D_{1,r}}\) the ordered sweep at time
\(r\),}
\[ \chi_{G_t}(q)\;=\;\sum_{c_H,\,c_O}\;\operatorname{tr}_{\mathbb C\mathcal P(c_H,c_O)}\Bigl(T^{(t-1)}\cdots T^{(1)}T^{(0)}\Bigr) \qquad\text{for every integer }q\ge0. \tag{T0} \]

{Proof.} Both sides count structures, and we exhibit a weight-preserving
bijection. Expand the trace: it is the sum, over pole colors and over
\(t\) proper path colorings \(\mathbf x^{(0)},\dots,\mathbf x^{(t-1)}\)
(indices mod \(t\)), of the product over all \((i,r)\) of
\(W_{D_{i,r}}\) evaluated at
\[ l=x_{i-1}^{(r+1)},\qquad x=x_i^{(r)},\qquad r=x_{i+1}^{(r)},\qquad X=x_i^{(r+1)} . \tag{T0a} \]
Three bookkeeping points make (T0a) correct. First, the ordering: when
\(T_i\) acts within sweep \(r\), sites \(1,\dots,i-1\) already carry
their time-\((r+1)\) values and sites \(i+1,\dots,n\) still carry
time-\(r\) values, so the left neighbor entering \(W\) is
\(x_{i-1}^{(r+1)}\) and the right neighbor is \(x_{i+1}^{(r)}\) ---
exactly the terminals of the face \(Q_{i,r}\) in (Q1). Second, the
poles: \(i-1=0\) or \(i+1=n+1\) produce the frozen colors \(c_H,c_O\),
matching the conventions \(v_{0,r}=H\), \(v_{n+1,r}=O\); the poles are
parameters of the trace, never dynamical sites, and their degree in
\(G_t\) plays no role. Third, the cyclic closure of the trace identifies
\(\mathbf x^{(t)}=\mathbf x^{(0)}\).

A proper coloring of \(G_t\) consists of the pole colors, the colors
\(x_i^{(r)}\) of all meridian vertices, and one internal extension in
each disk. Intermediate proper path states enforce exactly the
boundary-edge constraints of the quadrilaterals in (Q1). A shared
boundary edge can occur in two factors, but its inequality indicator is
idempotent and contributes no multiplicity. For an \(F_1\) cell, the
diagonal \(xX\) is an internal edge of that cell and is imposed by its
kernel \(W_{F_1}\). Internal vertices are disjoint between cells, so
their extension counts multiply. Thus the expanded trace and the
coloring count agree term by term. \(\square\)

By Lemma 2.1, \(T_i^{D}=uI+D_0U_i+(\mathsf B-D_0)(U_i^2-aU_i)\) on
\(\mathbb C\mathcal P(c_H,c_O)\), where \(U_i\) is the proper-recoloring
sum at site \(i\). The identity (T0) is therefore an identity for words
in \(I\) and the \(U_i\); what remains is to recognize those words
inside the Temperley--Lieb algebra and to evaluate the trace.

\hypertarget{the-coloring-representation-of-mathrmtl_n}{%
\subsection{\texorpdfstring{The coloring representation of
\(\mathrm{TL}_N\)}{The coloring representation of \textbackslash mathrm\{TL\}\_N}}\label{the-coloring-representation-of-mathrmtl_n}}

\textbf{Lemma 4.2 (coloring representation).} The operators in (T2a)
define a Temperley--Lieb representation; the following computation
proves the statement.

Let \(N=2n+2\) and let \(W_{\mathrm{col}}=(\mathbb C^q)^{\otimes(n+2)}\)
be the space of \emph{all} colorings of \(H,x_1,\dots,x_n,O\), proper or
not; note
\(\sum_{c_H,c_O}\operatorname{tr}_{\mathbb C\mathcal P(c_H,c_O)}(\cdot)\)
is the trace over \(W_{\mathrm{col}}\) of the same word multiplied by
the properness projector defined below. Set \(\delta=\sqrt q\) and
define, for \(0\le k\le n\) and \(1\le i\le n\),

\[ \rho(e_{2k})\;=\;\delta\,\Pi_k,\qquad \rho(e_{2i-1})\;=\;\delta^{-1}A_i, \tag{T2a} \]

where \(\Pi_k\) projects onto colorings with the \(k\)-th path edge
monochromatic (\(w_k=w_{k+1}\) in the vertex list
\(w_0=H,\dots,w_{n+1}=O\)) and \(A_i\) replaces the color at \(x_i\) by
a free sum over all \(q\) colors. This is the classical Potts
representation; we verify the relations rather than cite them.
\(\rho(e_{2k})^2=\delta^2\Pi_k=\delta\rho(e_{2k})\) and
\(\rho(e_{2i-1})^2=\delta^{-2}qA_i=\delta\rho(e_{2i-1})\) since
\(A_i^2=qA_i\). Nonadjacent generators commute: the even generators are
diagonal (even when they share a vertex), distinct odd generators act on
different sites, and a nonadjacent even--odd pair acts on disjoint
factors. For the braid-like relations, only adjacent pairs
\((e_{2i-1},e_{2i})\) and \((e_{2i-1},e_{2i-2})\) occur: on the factor
\(\mathbb C^q\otimes\mathbb C^q\) of \((x_i,w_{i+1})\),
\[ A_i\Pi_iA_i=A_i,\qquad \Pi_iA_i\Pi_i=\Pi_i, \] by direct computation
--- the first because the equalizer retains, out of the free sum over
the color at \(x_i\), exactly the one term matching the fixed neighbor,
after which the second resummation reproduces the original free sum; the
second because equalize--resum--equalize returns the equalizer. After
the normalizations in (T2a) these are precisely the relations
\(e_ae_be_a=e_a\):
\[ \bigl(\delta^{-1}A_i\bigr)\bigl(\delta\Pi\bigr)\bigl(\delta^{-1}A_i\bigr)=\delta^{-1}\,A_i\Pi A_i=\delta^{-1}A_i,\qquad \bigl(\delta\Pi\bigr)\bigl(\delta^{-1}A_i\bigr)\bigl(\delta\Pi\bigr)=\delta\,\Pi A_i\Pi=\delta\Pi. \]
So \(\rho\) is a representation of \(\mathrm{TL}_N(\delta)\) for every
\(\delta\ne0\) with \(\delta^2=q\).

\begin{figure}
\centering
\includegraphics[width=6.25in,height=\textheight]{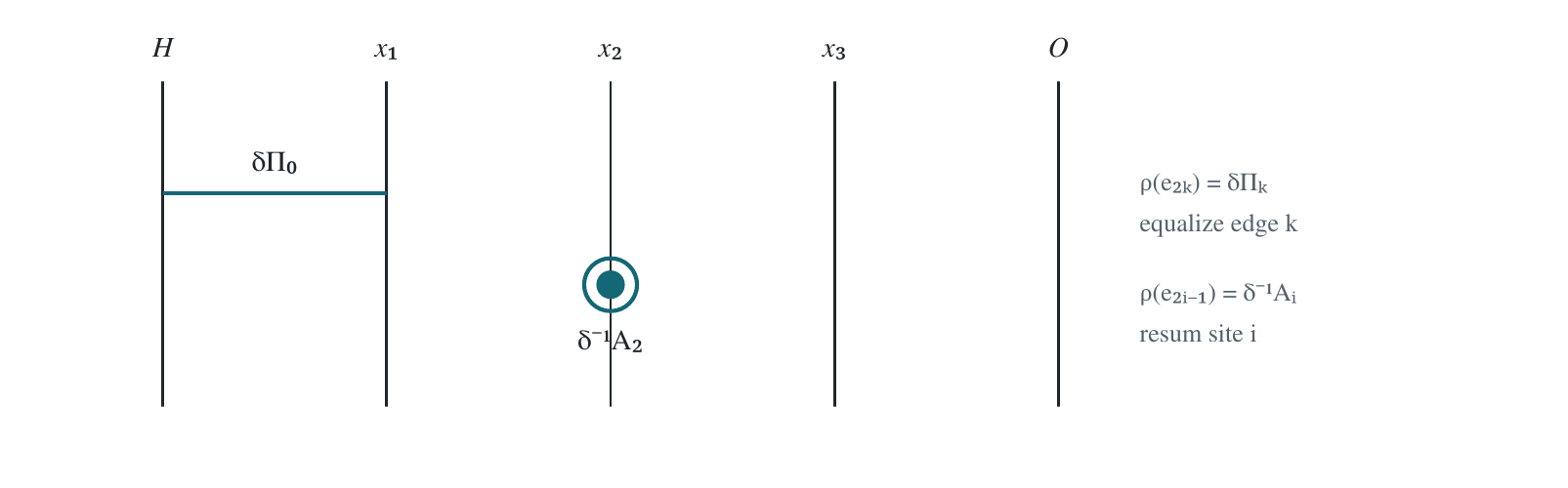}
\caption{The coloring representation of Lemma 4.2 on the vertex
worldlines of \(H,x_1,\dots,x_n,O\) (shown \(n=3\), \(N=8\)): even
generators are edge equalizers \(\delta\Pi_k\), odd generators are free
resummations \(\delta^{-1}A_i\) of a single site. The poles carry no odd
generator --- they have no recoloring move --- and the fused idempotent
\(F=\prod(I-e_{2k}/\delta)\) is exactly the properness projector.}
\end{figure}

Define \(F=\prod_{k=0}^{n}(I-e_{2k}/\delta)\) . Explicitly,

\[ f_k=I-e_{2k}/\delta,\qquad F=\prod_{k=0}^{n}f_k,\qquad F^2=F.\tag{T1}\]

Then \(\rho(F)=\prod_k(I-\Pi_k)\) is exactly the projector onto proper
path colorings: the properness of the state space in Lemma 4.1 is not an
analogy but the image of the fused idempotent. Moreover

\[ \rho\bigl(\delta\,Fe_{2i-1}F\bigr)\;=\;\rho(F)\,A_i\,\rho(F)\;=\;U_i, \tag{T2b} \]

since summing freely over the color at \(x_i\) and projecting to
properness on both sides is, by definition, the proper-recoloring
operator. Thus the sweep word of Lemma 4.1, conjugated into
\(W_{\mathrm{col}}\), is the \(\rho\)-image of the abstract
Temperley--Lieb element

\[ \mathsf U_i=\delta Fe_{2i-1}F,\qquad
\mathsf T_{i,r}=u_{i,r}F+D_{0,i,r}\mathsf U_i+(\mathsf B_{i,r}-D_{0,i,r})(\mathsf U_i^2-a\mathsf U_i),\]
\[ X_t=(\mathsf T_{n,t-1}\cdots\mathsf T_{1,t-1})\cdots(\mathsf T_{n,0}\cdots\mathsf T_{1,0}). \tag{T2c}\]

The unit of the corner algebra is \(F\); the latest sweep is written on
the left.

One sometimes streamlines (T2c) using the local cubic relation (D2)
inside the corner algebra \(F\,\mathrm{TL}_N F\); two different roles of
(D2) must be kept apart. For the trace identity of Theorem 4.4, no
abstract form of (D2) is needed: the identity is applied to the explicit
element (T2c) as it stands, and (D2) on the coloring space is Section
2\textquotesingle s theorem. For the sector norm bounds of Section 6,
however, the cubic constraint on the operators (T3) acting on the fused
modules \emph{is} needed, and the coloring-space proof does not supply
it there; Lemma 4.6 below proves it in the fused algebra itself. The
sector formulas (T3) used in Sections 5--6 compute the restriction of
this same element to the fused modules, diagram by diagram.

\hypertarget{evaluation-of-the-coloring-trace}{%
\subsection{Evaluation of the coloring
trace}\label{evaluation-of-the-coloring-trace}}

{Lemma 4.3 (trace evaluation).} \emph{For every
\(Y\in\mathrm{TL}_N(\delta)\) with \(\delta^2=q\), \(q\) a positive
integer,}
\[ \operatorname{tr}_{W_{\mathrm{col}}}\rho(Y)\;=\;q\cdot\operatorname{Tr}_{\mathrm{TL}}(Y), \tag{T2d} \]
\emph{where \(\operatorname{Tr}_{\mathrm{TL}}\) assigns to each diagram
\(\delta^{\#\text{loops of its trace closure}}\) and extends linearly.}

{Proof.} We argue by induction on the number of strands, through the
Markov property that characterizes the diagram trace; the identity to be
proved is the top case of a one-parameter family. For every \(M\ge1\)
define a coloring representation \(\rho_M\): for even \(M=2m+2\), on the
\(m+2\) vertices \(w_0,\dots,w_{m+1}\) exactly as in (T2a), both end
vertices carrying no resummation generator; for odd \(M=2m+1\), on the
\(m+1\) vertices \(w_0,\dots,w_m\), with edge equalizers \(\delta\Pi_k\)
at the even indices and resummations \(\delta^{-1}A_i\) at \emph{all} of
\(w_1,\dots,w_m\), so that the last generator \(e_{M-2}=\delta^{-1}A_m\)
acts at the right end; the defining relations are verified by the same
two-factor computations as in Section 4.2, applied at the right end. Set
\(\varphi_M=\operatorname{tr}\circ\rho_M\) and
\[ c_M=\begin{cases}q,&M\text{ even},\\ \delta,&M\text{ odd}.\end{cases} \]
We prove \(\varphi_M=c_M\operatorname{Tr}_{\mathrm{TL}_M}\) for every
\(M\); the case \(M=N\) is (T2d).

\textbf{A spanning lemma.} The standard basic-construction span is
\(\mathrm{TL}_M=\mathrm{TL}_{M-1}+\mathrm{TL}_{M-1}e_{M-2}\mathrm{TL}_{M-1}\).
For completeness, induct on \(M\). In a word with two occurrences of the
top generator \(e=e_{M-2}\), expand the intervening element of
\(\mathrm{TL}_{M-1}\), by induction, as a sum of elements \(x\) and
\(y e_{M-3}z\) with \(x,y,z\in\mathrm{TL}_{M-2}\). These commute with
\(e\), so \(exe=\delta xe\) and \(e y e_{M-3}z e=yz e\). Repeating
reduces the number of top generators to at most one. The base cases are
immediate. This is also the usual uniqueness argument for the Markov
trace {[}LK, Theorem 5.2{]}.

\textbf{Three shared properties.} Both functionals \(\varphi_M/c_M\) and
\(\operatorname{Tr}_{\mathrm{TL}_M}\) satisfy: (0) the value
\(\delta^M\) on the identity --- for \(\varphi\), count vertices:
\(q^{m+2}=q\,\delta^{M}\) in the even case and
\(q^{m+1}=\delta\,\delta^{M}\) in the odd case; (1) the trace property
\(F(XY)=F(YX)\); and (2) the two Markov relations, for
\(x,z\in\mathrm{TL}_{M-1}\):
\[ F_M(x)=\delta\,F_{M-1}(x),\qquad F_M\bigl(z\,e_{M-2}\bigr)=F_{M-1}(z). \tag{T2e} \]
For the diagram trace, (T2e) is the standard isotopy pair: restricting
to \(\mathrm{TL}_{M-1}\) leaves the last strand as a free closed loop,
and closing \(z\,e_{M-2}\) merges that loop through the cap into the
closure of \(z\). For the coloring side, (T2e) is a direct computation,
in which the parity constants earn their keep. Even \(M\): a word \(x\)
avoiding \(e_{M-2}=\delta\Pi_{m}\) never touches the end vertex
\(w_{m+1}\), so \(\varphi_M(x)=q\,\varphi_{M-1}(x)\), and
\(q/c_M=\delta/c_{M-1}\); while \(\varphi_M(z\,\delta\Pi_m)\) ties
\(w_{m+1}\) to \(w_m\) and sums it, giving \(\delta\,\varphi_{M-1}(z)\),
and \(\delta/c_M=1/c_{M-1}\). Odd \(M\): a word \(x\) avoiding
\(e_{M-2}=\delta^{-1}A_m\) is represented by exactly the same operators
in \(\rho_M\) and in \(\rho_{M-1}\), on the same space, so
\(\varphi_M(x)=\varphi_{M-1}(x)\), and \(1/c_M=\delta/c_{M-1}\) since
\(q=\delta^2\). For the second relation no tensor factorization is
available --- the equalizer \(e_{M-3}=\delta\Pi_{m-1}\) acts on \(w_m\)
--- and the correct observation is weaker and suffices: every generator
available to \(z\) is \emph{diagonal in the color of \(w_m\)}, since
equalizers change no colors and the resummations in \(z\) act at
\(w_1,\dots,w_{m-1}\) only. Writing the matrix element of \(A_m\) as
\(\langle c',b|A_m|c'',b'\rangle=\mathbf 1_{c'=c''}\) (the new color
\(b\) is free, the old \(b'\) is forgotten),
\[ \operatorname{tr}\bigl(\rho_M(z)A_m\bigr)=\sum_{c',\,b,\,b'}\langle c',b\,|\,\rho_M(z)\,|\,c',b'\rangle =\sum_{c',\,b}\langle c',b\,|\,\rho_M(z)\,|\,c',b\rangle=\operatorname{tr}\rho_M(z), \]
the middle step because diagonality in the last coordinate kills every
term with \(b'\ne b\). Hence
\(\varphi_M(z\,\delta^{-1}A_m)=\delta^{-1}\varphi_M(z)=\delta^{-1}\varphi_{M-1}(z)\),
and \(\delta^{-1}/c_M=1/c_{M-1}\). All eight ratios match.

\textbf{Induction.} Suppose
\(\varphi_{M-1}=c_{M-1}\operatorname{Tr}_{\mathrm{TL}_{M-1}}\). By the
spanning lemma, every element of \(\mathrm{TL}_M\) is a combination of
elements \(x\) and \(X e_{M-2}Y\) with \(x,X,Y\in\mathrm{TL}_{M-1}\).
Property (1) gives \(F_M(Xe_{M-2}Y)=F_M(YXe_{M-2})\), and then (T2e)
determines both kinds of values from \(F_{M-1}\). Two functionals
satisfying (0)--(2) with the same \(F_{M-1}\) therefore agree on
\(\mathrm{TL}_M\). The base \(M=1\) is the one-dimensional algebra.
{\(\square\)}

Lemma 4.3 settles, in particular, that the trace functional appearing in
the chromatic correspondence and the standard quantum trace decomposed
in (T4) are \emph{the same functional}, with the overall factor \(q\)
accounting for the two pole regions of the sphere; no separate matching
of normalizations between {[}FK1{]} and the Markov-trace literature
needs to be assumed.

\hypertarget{fused-bases-without-singular-projectors}{%
\subsection{Fused bases without singular
projectors}\label{fused-bases-without-singular-projectors}}

For \(0\le j\le n+1\), let \(W_j=W_{N,2j}\) be the standard link-state
module with \(2j\) defects and \(p=n+1-j\) caps: link states are
noncrossing caps and defects extending outward, caps may not enclose a
defect, joining two defects gives zero, and a closed loop contributes
\(\delta\). Keep exactly those states \(b\) containing none of the caps
\((2k,2k+1)\). The vectors \(Fb\) form a basis of \(V_j=FW_j\): a
discarded state is annihilated by the relevant factor of \(F\); and in
the expansion of a retained \(Fb\), the coefficient of \(b\) itself is
\(1\) while every correction term carries a forbidden cap, so the
restriction \(R_0\) to retained coordinates is exactly the identity ---
independence and spanning at every \(\delta\ne0\). The local operators
become

\[ U_i=\delta\,R_0\,e_{2i-1}\Bigl(I-\frac{e_{2i-2}}{\delta}\Bigr)\Bigl(I-\frac{e_{2i}}{\delta}\Bigr)\big|b\rangle , \tag{T3} \]

with four diagrammatic terms and Laurent-polynomial entries. The top
fused dimensions are \(\dim V_{n+1}=1\), \(\dim V_n=n\),
\(\dim V_{n-1}=n(n+1)/2\). The certificate code constructs these bases
literally rather than assuming their dimensions.

The decomposition of the diagram trace into sector traces is classical;
because it carries real weight in this paper, we prove it rather than
cite it (compare {[}J{]}).

{Proposition 4.5 (sector decomposition of the diagram trace).} \emph{For
every \(Y\in\mathrm{TL}_N(\delta)\),}
\[ \operatorname{Tr}_{\mathrm{TL}}(Y)\;=\;\sum_{j=0}^{n+1}\,[2j+1]_\delta\;\operatorname{tr}_{W_j}(Y). \tag{T4} \]

{Proof.} Fix a word \(Y\) in the generators with Laurent-polynomial
coefficients; both sides of (T4) are then Laurent polynomials in
\(\delta\), so it suffices to prove the identity for every real
\(\delta\gt2\), where three structural facts hold. First, every Gram
determinant (G1) is a product of quotients \([m]_\delta/[r]_\delta\)
with positive numerators and denominators --- for \(\delta\gt2\) write
\(\delta=2\cosh\theta\), \([m]_\delta=\sinh(m\theta)/\sinh\theta\gt0\)
--- so by the inertia-constancy lemma (Lemma 6.0 below) each Gram form
is positive definite, and each standard module \(W_j\) is simple by the
rank-one spanning argument of Section 6.1. Second, their direct sum is
faithful. To see this, filter the diagram basis by cap number and take
the smallest cap number \(p\) occurring in a nonzero element. Diagrams
with exactly \(p\) caps are pairs of upper and lower \(p\)-cap link
states; let \(A\ne0\) be their coefficient matrix. On the standard
module with \(p\) caps the diagrams with more caps vanish, and this
leading layer acts by \(A G_{N,p}\), which is nonzero since the Gram
matrix is invertible. Thus no nonzero element annihilates all standard
modules. Third, the dimensions satisfy the classical Catalan convolution
\(\sum_j(\dim W_j)^2=\dim\mathrm{TL}_N\), so the faithful map
\(\mathrm{TL}_N\to\bigoplus_j\operatorname{End}(W_j)\) is an
isomorphism, and every trace functional on \(\mathrm{TL}_N\) is a unique
linear combination \(\sum_jc_j\operatorname{tr}_{W_j}\).

It remains to compute the coefficients for the diagram trace. Evaluate
both sides on the disjoint-cap elements \(E_p=e_0e_2\cdots e_{2p-2}\)
(empty product for \(p=0\)), \(0\le p\le n+1\). The trace closure of
\(E_p\) has \(N-p\) loops, so the left side is \(\delta^{\,N-p}\). On
\(W_j\), the product \(E_p\) maps a link state \(b\) to a multiple of a
state containing the caps \((0,1),\dots,(2p-2,2p-1)\); the diagonal
entry therefore vanishes unless \(b\) contains all \(p\) caps, in which
case it equals \(\delta^{\,p}\), so
\(\operatorname{tr}_{W_j}(E_p)=\delta^{\,p}\dim W_{N-2p,\,2j}\), the
count of noncrossing completions on the remaining \(N-2p\) sites. Since
\(\dim W_{M,2j}=0\) for \(2j\gt M\) and the extreme dimension is \(1\),
the matrix \(\bigl(\operatorname{tr}_{W_j}(E_p)\bigr)_{p,j}\) is
triangular with nonzero diagonal: the coefficients \(c_j\) are
determined by these evaluations alone. Finally, \(c_j=[2j+1]_\delta\)
satisfies all of them at once, by the following identity with
\(M=N-2p\): writing \(\delta=z+z^{-1}\), so that
\([2j+1]_\delta=\sum_{r=-j}^{j}z^{2r}\), and
\(\dim W_{M,2j}=\binom{M}{M/2-j}-\binom{M}{M/2-j-1}\),

\[ \begin{aligned}\sum_{j\ge0}[2j+1]_\delta\dim W_{M,2j} &=\sum_{r\in\mathbb Z}z^{2r}\!\!\sum_{j\ge|r|}\!\Bigl(\tbinom{M}{M/2-j}-\tbinom{M}{M/2-j-1}\Bigr) \\ &=\sum_{r\in\mathbb Z}z^{2r}\tbinom{M}{M/2-|r|} =\bigl(z+z^{-1}\bigr)^{M}=\delta^{M}.\end{aligned}\tag{T4a} \]

the inner sum telescoping and the last step being the binomial theorem.
Hence
\[\delta^{\,p}\sum_jc_j\dim W_{N-2p,2j}=\delta^{\,p}\,\delta^{\,N-2p}=\delta^{\,N-p}.\]
for every \(p\), as required. {\(\square\)}

Both sides of (T4), evaluated on a fixed word, being Laurent polynomials
in \(\delta\), the identity holds at every parameter, in particular at
every Beraha value; this use of generic semisimplicity, like that in
Lemma 4.6, enters as a polynomial identity, never as a structural claim
at the target. For \(Y=FYF\) one has
\(W_j=\operatorname{im}F\oplus\ker F\) with \(Y\) vanishing on the
second summand, so for positive powers the trace in (T4) equals the
trace on \(V_j\).

{Lemma 4.6 (the cubic relation in the fused algebra).} \emph{In
\(F\,\mathrm{TL}_N(\delta)F\), the elements \(u_i=\delta Fe_{2i-1}F\)
satisfy}
\[ u_i\,\bigl(u_i-aF\bigr)\,\bigl(u_i-(a+1)F\bigr)\;=\;0,\qquad a=\delta^2-2, \tag{T3a} \]
\emph{for every \(\delta\ne0\); consequently the sector operators (T3)
satisfy the same cubic on every \(V_j\). For \(q=\delta^2\gt2\),
\(E_i=U_i(U_i-aI)/(a+1)\) is an idempotent there.}

{Proof.} Locality first. Write \(F=F_{\mathrm{loc}}F_{\mathrm{rest}}\)
with \(F_{\mathrm{loc}}=f_{i-1}f_i\); every factor of
\(F_{\mathrm{rest}}\) commutes with \(e_{2i-1}\), with
\(F_{\mathrm{loc}}\), and with each other, and is idempotent, so
\(u_i=F_{\mathrm{rest}}\,v\,F_{\mathrm{rest}}\) with
\(v=\delta F_{\mathrm{loc}}e_{2i-1}F_{\mathrm{loc}}\), and powers obey
\(u_i^k=F_{\mathrm{rest}}\,v^k F_{\mathrm{rest}}\). It therefore
suffices to prove the cubic for \(v\), which is the image of the
corresponding element of \(\mathrm{TL}_4\) under the strand embedding
onto positions \(2i-2,\dots,2i+1\). In \(\mathrm{TL}_4\) at a fixed real
\(\delta\gt2\), the sum of the standard modules is faithful (first
paragraph of the proof of Proposition 4.5), so it suffices to verify the
cubic on each fused sector \(V_j\), \(j=0,1,2\). All three are
one-dimensional, and a direct diagram computation evaluates \(v\) on
them as the scalars
\[ v\big|_{V_0}=a+1,\qquad v\big|_{V_1}=a,\qquad v\big|_{V_2}=0, \] the
three roots of the cubic. The three evaluations, written out: on
\(V_2\), the all-defect state is annihilated by every generator, so
\(v=0\) there. On \(V_1\), the retained state \(c\) has the single cap
\((1,2)\); the two forbidden-position neighbors are \(c^-\) with cap
\((0,1)\) and \(c^+\) with cap \((2,3)\), and
\(F_{\mathrm{loc}}c=c-\delta^{-1}c^--\delta^{-1}c^+\). Applying \(e_1\):
the cap \((1,2)\) of \(c\) gives \(\delta c\), and each neighbor rewires
back to \(c\) with coefficient \(1\), so
\(e_1F_{\mathrm{loc}}c=(\delta-2\delta^{-1})c=(a/\delta)c\) and
\(v\,c=\delta F_{\mathrm{loc}}(a/\delta)c=a\,F_{\mathrm{loc}}c\): the
value is \(a\). On \(V_0\), the retained state \(b\) has the nested caps
\((0,3),(1,2)\); both \(e_0\) and \(e_2\) send it to the same forbidden
state \(f\) with caps \((0,1),(2,3)\), and
\(F_{\mathrm{loc}}b=b-\delta^{-1}f\) after the cross term is absorbed by
\(e_2f=\delta f\). Applying \(e_1\): \(e_1b=\delta b\) and \(e_1f=b\),
so \(e_1F_{\mathrm{loc}}b=(\delta-\delta^{-1})b=((a+1)/\delta)\,b\) and
the value is \(a+1\). (The supplied symbolic trace audit verifies all
three evaluations identically in \(\delta\).) Hence (T3a) holds in
\(\mathrm{TL}_4\) for every real \(\delta\gt2\); since each side has
Laurent-polynomial coefficients in the diagram basis, it holds for every
\(\delta\ne0\), and its image under the embedding and the
\(F_{\mathrm{rest}}\)-dressing gives (T3a) in \(\mathrm{TL}_N\).
Restriction to the invariant subspaces \(V_j\) is a homomorphism,
proving the sector statement; for \(q\gt2\) its three distinct roots
give \(E_i^2=E_i\) by polynomial interpolation. {\(\square\)}

\hypertarget{the-boxed-identity}{%
\subsection{The boxed identity}\label{the-boxed-identity}}

Define \([0]_\delta=0\), \([1]_\delta=1\),
\([k+1]_\delta=\delta[k]_\delta-[k-1]_\delta\); for
\(\delta=2\cos\theta\), \([k]_\delta=\sin(k\theta)/\sin\theta\).
Whenever the fixed disk has \(u\lt0\), define the following
normalization (the strict filter inequalities are not required for this
identity): \(g=-u\), \(r_v=v/g\), \(r_e=e/g\), \(\eta=(1+r_v)/a\),
\(E_i=U_i(U_i-aI)/(a+1)\),
\(S_i=-T_i/g=I-\eta U_i-(1+r_e-\eta(a+1))E_i\), and
\(P_j=(S_n\cdots S_1)|_{V_j}\). The identity in each fused module is
denoted by \(I\). In particular,
\[ S_i=I-\eta U_i-\bigl(1+r_e-\eta(a+1)\bigr)E_i,\qquad P_j=(S_n\cdots S_1)|_{V_j}. \tag{T5} \]

{Theorem 4.4 (complete fused trace).} \emph{For the uniform closure
\(G_t=G_t(D,n)\) of a disk with eigenvalue scale \(g\),}
\[ \boxed{\;\chi_{G_t}(q)\;=\;(-g(q))^{\,nt}\,q\sum_{j=0}^{n+1}\,[2j+1]_{\sqrt q}\;\operatorname{tr}\!\bigl(P_j(q)^{\,t}\bigr),\;} \tag{T6} \]
\emph{and more generally, for arbitrary cells \(D_{i,r}\), the same
identity holds with \((-g)^{nt}\operatorname{tr}(P_j^t)\) replaced by
the trace of the ordered product of the unnormalized sector transfers.
The unnormalized identity holds identically as a Laurent-polynomial
identity in \(\delta\), and its sum is the polynomial
\(\chi_{G_t}(\delta^2)\). The normalized formula is used on
neighborhoods where \(q\gt2\) and \(g\ne0\).}

{Proof.} For each positive integer \(q\), Lemmas 4.1--4.3 identify
\(\chi_{G_t}(q)\) with \(q\operatorname{Tr}_{\mathrm{TL}}X_t\). For a
fixed finite disk and closure, both expressions are Laurent polynomials
in \(\delta\): the disk coefficients are polynomials in \(q=\delta^2\),
and the only denominators in the unnormalized TL element are powers of
\(\delta\). Their difference vanishes at infinitely many values
\(\delta=\sqrt q\), and therefore is identically zero after
multiplication by a sufficiently large power of \(\delta\). Apply (T4)
and restrict to the fused modules. In the uniform case each local
transfer is \(-gS_i\), giving (T6) wherever the normalized expressions
are defined. In particular all of these expressions are real analytic in
a fixed neighborhood of each target considered below. For \(t=2m\),
\((-g)^{nt}q\) is positive there. {\(\square\)}

The identity transfer gives
\(q\operatorname{Tr}_{\mathrm{TL}}F=q(q-1)^{n+1}\), the number of proper
path colorings. Appendix B checks the conventions against nine
independently computed graph polynomials, including nonuniform rotated
and serially composed cells. These finite checks supplement the general
proof above; they are not an independent review of it.

\hypertarget{s5}{%
\section{A carrier with one expanding eigenvalue}\label{s5}}

The sector that carries the mechanism is the one-cap sector \(j=n\).
Using as basis the \(n\) fused caps between consecutive pairs, in their
spatial order, formula (T3) collapses to a strikingly simple form:

\[ U_i=\mathbf e_i\,\mathbf e_i^{T}H, \qquad H \;=\; aI-\operatorname{Adj}(P_n), \tag{H1} \]

where \(\operatorname{Adj}(P_n)\) is the adjacency matrix of the path on
\(n\) vertices: the diagonal coefficient is \(\delta^2-2=a\), the two
neighbors contribute \(-1\), and everything else vanishes. Written out:
the retained states are \(c_1,\dots,c_n\), where \(c_k\) has its single
cap at the odd pair \((2k-1,2k)\), and
\(Fc_k=c_k-\delta^{-1}c_k^--\delta^{-1}c_k^+\), with \(c_k^{\pm}\) the
forbidden even-position shifts of the cap. Then \(e_{2i-1}\) returns
\(\delta c_i\) on \(c_i\), returns \(c_i\) with coefficient \(1\) on
each of \(c_{i-1}^+\) and \(c_{i+1}^-\), and annihilates every other
term in the standard module (joining two defects raises the cap count);
collecting coefficients, \(U_i c_k=\delta Fe_{2i-1}Fc_k\) has diagonal
entry
\(\delta\cdot(\delta-2\delta^{-1})\cdot\langle\text{retained}\rangle=a\)
at \(k=i\), entry \(-1\) at \(k=i\pm1\) through the \(-\delta^{-1}\)
components, and zero beyond. In particular \(U_i^2=aU_i\) on this
sector, so the projector term \(E_i\) vanishes and the sweep is built
purely from rank-one updates
\(S_i=I-\eta\,\mathbf e_i\mathbf e_i^{T}H\).

{Lemma 5.1.} \emph{Suppose \(a\gt0\), \(H\) is nonsingular with inertia
\((n-1,1)\), and \(0\lt r_v\lt1\). Then every word \(P\) in the
operators \(S_i=I-\eta\,\mathbf e_i\mathbf e_i^{T}H\) that uses every
site has exactly one eigenvalue outside the unit circle. That eigenvalue
is positive, real, and algebraically simple; every other eigenvalue has
modulus strictly less than \(1\).}

{Proof.} The engine is a one-line congruence. Since
\(a\eta=1+r_v\in(1,2)\),
\[ H-S_i^{T}HS_i \;=\; \eta\,(2-a\eta)\,H\mathbf e_i\mathbf e_i^{T}H \;\succeq\;0 , \tag{H2} \]
and telescoping along the word gives \(Q=H-P^{T}HP\succeq0\). Suppose
\(z^{T}Qz=0\). Writing the sweep word as \(P=S_{i_L}\cdots S_{i_1}\)
with partial products \(P_r=S_{i_r}\cdots S_{i_1}\) and \(P_0=I\), the
telescoping expresses \(z^TQz\) as the sum over \(r\) of the increments
\((P_{r-1}z)^{T}\bigl(H-S_{i_r}^{T}HS_{i_r}\bigr)(P_{r-1}z)\), each
nonnegative by (H2), so each vanishes. By (H2) the \(r\)-th increment is
a positive multiple of \(\bigl(\mathbf e_{i_r}^{T}HP_{r-1}z\bigr)^2\),
so \(\mathbf e_{i_r}^{T}HP_{r-1}z=0\); and since the update \(S_i\)
changes a vector only through the scalar \(\mathbf e_i^{T}H(\cdot)\),
this forces \(P_rz=P_{r-1}z\). By induction \(P_rz=z\) and
\(\mathbf e_{i_r}^{T}Hz=0\) for every \(r\); every site occurs in the
word, so \(Hz=0\), contradicting the nonsingularity of \(H\). Hence
\(Q\gt0\).

There is no unit-modulus eigenvalue: a complex eigenvector \(z\) with
\(|\lambda|=1\) would give \(z^{*}Qz=(1-|\lambda|^2)z^{*}Hz=0\). Each
\(S_i\) is invertible, with determinant \(-r_v\ne0\). On the stable
generalized eigenspace, telescoping forward,
\[ x^{T}Hx=\sum_{k\ge0}(P^{k}x)^{T}Q(P^{k}x)\gt0\quad(x\ne0), \] while
on the unstable generalized eigenspace, telescoping backward,
\[ x^{T}Hx=-\sum_{k\ge1}(P^{-k}x)^{T}Q(P^{-k}x)\lt0 . \] The two real
invariant subspaces exhaust the space, so the inertia \((n-1,1)\) of
\(H\) forces their dimensions to be \(n-1\) and \(1\): exactly one
unstable eigenvalue, necessarily real and algebraically simple. For
positivity, consider for \(0\le s\le1\) the path
\(I-s\eta\,\mathbf e_i\mathbf e_i^{T}H\); by (H2) with \(\eta\) replaced
by \(s\eta\), each map along the path sends every \(H\)-negative vector
to an \(H\)-negative vector, so the path preserves each component of the
negative cone, hence so does the product, and its unstable eigenvalue is
positive. {\(\square\)}

At \(q_*=B_{2n+1}\) or \(B_{2n}\) with \(n\ge3\), the eigenvalues of
\(H\) are \(a-2\cos\!\frac{k\pi}{n+1}\), \(1\le k\le n\); the first is
negative and the rest positive, so Lemma 5.1 applies whenever
\(0\lt r_v\lt1\).

For the ascending sweep one can say more. Writing \(H=aI+L_-+L_+\) with
strictly lower and upper triangular parts, sequential updating gives
\(P=(I+\eta L_-)^{-1}(-r_vI-\eta L_+)\), whose characteristic polynomial
obeys the three-term recurrence

\[ p_0=1,\quad p_1=\lambda+r_v,\quad p_k=(\lambda+r_v)\,p_{k-1}-\eta^{2}\lambda\,p_{k-2}, \tag{H3} \]

and the expanding eigenvalue has the closed form, with
\(c=\cos(\pi/(n+1))\),

\[ \lambda_* \;=\; \Bigl(\eta c+\sqrt{\eta^{2}c^{2}-r_v}\Bigr)^{2}. \tag{H4} \]

To derive (H4) from (H3): the linear recurrence has characteristic roots
\(\zeta_\pm\) with \(\zeta_++\zeta_-=\lambda+r_v\) and
\(\zeta_+\zeta_-=\eta^2\lambda\), and
\(p_n=(\zeta_+^{\,n+1}-\zeta_-^{\,n+1})/(\zeta_+-\zeta_-)\); this
vanishes exactly when \(\zeta_+/\zeta_-\) is a nontrivial \((n{+}1)\)-st
root of unity, i.e. when
\(\lambda+r_v=2\eta\sqrt\lambda\,\cos\tfrac{k\pi}{n+1}\) for some
\(1\le k\le n\). Each value of \(k\) gives a quadratic in
\(\sqrt\lambda\); the largest root over all \(k\) occurs at \(k=1\) and
is (H4). The inertia range used here satisfies \(a\lt2c\): at
\(q_*=B_\ell\) one has \(a=2\cos(2\pi/\ell)\) and \(c=\cos(\pi/(n+1))\),
and \(2\pi/\ell\gt\pi/(n+1)\) amounts to \(\ell\lt2n+2\), true for both
\(\ell=2n+1\) and \(\ell=2n\). Consequently the quadratic
\((\lambda+r_v)^2-4\eta^2c^2\lambda\) is negative at \(\lambda=1\); its
larger root exceeds \(1\) and must be the unique expanding eigenvalue.
The same quadratic certifies lower bounds on \(\lambda_*\) in the exact
arithmetic of Section 9.

\hypertarget{s6}{%
\section{Gram forms and all competing sectors}\label{s6}}

\hypertarget{standard-input-and-the-positivity-test}{%
\subsection{Standard input and the positivity
test}\label{standard-input-and-the-positivity-test}}

Let \(G_{N,p}(\delta)\) be the diagram Gram matrix of the standard
module with \(p\) caps; the generators are self-adjoint for this form.
Two external inputs are imported from Ridout and Saint-Aubin {[}RS{]}
and used as black boxes, with all real inertia arguments supplied here.
The first is the determinant formula {[}RS, Theorem 4.7{]}:

\[ \det G_{N,p} \;=\; \prod_{r=1}^{p}\left(\frac{[\,N-2p+1+r\,]_\delta}{[\,r\,]_\delta}\right)^{d_{N,p-r}}, \qquad d_{N,s}=\binom Ns-\binom N{s-1}. \tag{G1} \]

Three conventions of {[}RS{]} are fixed once, because the dictionary of
Section 6.3 depends on them. In {[}RS{]}, an \((\mathsf n,p)\)-link
state has \(\mathsf n\) points and \(p\) links, hence \(\mathsf n-2p\)
defects: \emph{the second index counts links, not defects} {[}RS,
Section 2{]}. The standard module \(\mathcal V_{\mathsf n,p}\) is
spanned by these states, with the diagram action truncated to zero
whenever the number of links increases {[}RS, Equation (3.3){]}; its
radical \(\mathcal R_{\mathsf n,p}\) is by definition the radical of the
bilinear form, so it is nonzero exactly when the Gram determinant
vanishes. The pair \((\mathsf n,p)\) is critical if
\(\mathsf q^{2(\mathsf n-2p+1)}=1\), and two pairs \((\mathsf n,p')\),
\((\mathsf n,p)\) with \(0\lt|p'-p|\lt\ell_{\mathrm{RS}}\) form a
\emph{symmetric pair} if both are noncritical and they lie symmetrically
on either side of the single critical line between them {[}RS, Section
7{]}; here \(\ell_{\mathrm{RS}}\ge2\) is minimal with
\(\mathsf q^{2\ell_{\mathrm{RS}}}=1\).

{Imported Theorem ({[}RS, Theorem 7.2{]}).} \emph{The radical
\(\mathcal R_{\mathsf n,p}\) is zero or irreducible and, if the pair
\(\mathcal V_{\mathsf n,p}\) and \(\mathcal V_{\mathsf n,p'}\), with
\(p'\gt p\), is symmetric, then
\(\mathcal R_{\mathsf n,p'}\cong\mathcal L_{\mathsf n,p}\), the
irreducible quotient of its partner.}

This is the statement of {[}RS{]} with only notational adjustment; no
other consequence of {[}RS, Section 7{]} is used in this paper, and the
determinant formula (G1) is {[}RS, Theorem 4.7{]} verbatim.

One classical consequence of nondegeneracy is worth recalling: a diagram
with upper link state \(s\) and lower link state \(t\) acts on the
standard module as the rank-one map \(w\mapsto s\langle t,w\rangle\); if
the Gram form is nondegenerate these maps span all endomorphisms, so the
standard module is simple.

{Lemma 6.0 (inertia constancy).} \emph{A continuous family of real
symmetric matrices whose determinant does not vanish on an interval has
constant inertia there.}

{Proof.} Eigenvalues of a continuous symmetric family vary continuously;
if a signature changed, some eigenvalue would cross zero, making the
determinant vanish. {\(\square\)}

The diagonal entries of \(G_{N,p}\) are \(\delta^{p}\), and any nonzero
off-diagonal entry has degree at most \(p-1\); hence the real Gram
matrix is positive definite for all sufficiently large positive
\(\delta\), and on any interval where its determinant does not vanish
its inertia is constant. It is this constancy (Lemma 6.0) --- not the
sign of a determinant, which can hide an even number of negative
directions --- that establishes positivity throughout the arguments
below; every later appeal to "continuation from large \(\delta\)" is an
appeal to this lemma. The projector \(F\) is self-adjoint, so on a
positive module it induces a positive restricted form on \(V_j\); on a
semidefinite module,

\[ \operatorname{rad}\bigl(G|_{\operatorname{im}F}\bigr)=\operatorname{im}F\cap\operatorname{rad}G, \tag{G2} \]

since for \(v=Fv\), orthogonality to \(\operatorname{im}F\) already
implies orthogonality to every \(w\) through
\(\langle v,w\rangle=\langle v,Fw\rangle\). Finally, by Lemma 4.6 and
(T5), on any positive module or positive quotient,

\[ \|S_i\|\le M:=\max(1,|r_v|,|r_e|), \qquad \rho(P_j)\le M^{\,n}, \tag{G3} \]

all norms taken in the module\textquotesingle s own Gram form. Under the
filter condition (D8), \(M=1\): this is where Section 2 pays for Section
6.

\hypertarget{odd-indices}{%
\subsection{Odd indices}\label{odd-indices}}

{Proposition 6.1.} \emph{At \(q_*=B_{2n+1}\), \(n\ge3\), under
\(0\lt r_v,r_e\lt1\), every sector other than \(V_n\) has spectral
radius at most \(1\).}

{Proof.} Put \(\delta_*=2\cos(\pi/(2n+1))\). For \(j\le n-2\), with
\(p=n+1-j\), the largest numerator index in (G1) is
\(2j+1+p=n+j+2\le2n\); every quantum integer appearing is positive on
\([\delta_*,\infty)\) (the sine formula up to \(2\), the hyperbolic
formula beyond), so the Gram form is positive definite at the target and
(G3) applies with \(M=1\). For \(j=n-1\), \(p=2\), the determinant is
\[ \det G_{N,2}=[2n]^{\,N-1}\,\frac{[2n+1]}{[2]} , \tag{G4} \] nonzero
above \(\delta_*\) and with a simple zero at \(\delta_*\): the limiting
form is positive semidefinite with a one-dimensional radical, which is
invariant. On a one-dimensional TL module, writing \(e_i=x_i\), the
relations force \(x_i\in\{0,\delta_*\}\) and \(x_ix_{i+1}x_i=x_i\); a
nonzero \(x_i\) would force its neighbor to equal \(1/\delta_*\),
impossible since \(\delta_*^2\ne1\), and every generator has a neighbor.
So every \(e_i\) --- hence every \(U_i\) --- vanishes on the radical,
where each \(S_i\) acts as \(1\). The quotient is positive and obeys
(G3); a possible nonsplit extension does not change the union of
subspace and quotient spectra. The top sector \(j=n+1\) likewise has
\(U_i=0\) and \(P=1\). All sectors are now accounted for. {\(\square\)}

\hypertarget{even-indices-and-the-positive-quotient}{%
\subsection{Even indices and the positive
quotient}\label{even-indices-and-the-positive-quotient}}

The translation into {[}RS{]} indexing, with \(p\) counting links as
fixed in Section 6.1: our \(W_j\), spanned by states with \(2j\) defects
on \(N\) strands, is \(\mathcal V_{N,(N-2j)/2}\); in particular
\(W_n=\mathcal V_{N,1}\) and \(W_{n-1}=\mathcal V_{N,2}\). At
\(q_*=B_{2n}\) one has \(\mathsf q=e^{i\pi/2n}\), so
\(\ell_{\mathrm{RS}}=2n\); the two labels give \(N-2p+1=2n+1\) and
\(2n-1\), which are noncritical, straddle the single critical value
\(2n=\ell_{\mathrm{RS}}\), and lie symmetrically about it:
\((\mathcal V_{N,2},\mathcal V_{N,1})\) is a symmetric pair in the sense
fixed in Section 6.1. The imported theorem, in this instance and with
its direction displayed, reads
\(\operatorname{rad}\mathcal V_{N,2}\cong\mathcal L_{N,1}\): the radical
lives in the two-link module \(W_{n-1}\) and is isomorphic to the
irreducible head of the one-link module \(W_n\) --- this, and not the
reverse, is the direction consumed below at (G7). The excluded
parameters \(\mathsf q=\pm1\), i.e. \(\delta=\pm2\), do not arise at any
noninteger Beraha value.

The even case is where the algebra at the root of unity stops being an
obstacle and becomes the mechanism. The statement to be proved:

{Proposition 6.2.} \emph{At \(q_*=B_{2n}\), \(n\ge4\): (1) \(W_j\) is
positive for \(j\le n-3\); (2) \(W_{n-2}\) is positive semidefinite with
a one-dimensional radical on which all \(U_i\) vanish; (3) \(W_n\) is
simple and \(\operatorname{rad}W_{n-1}\cong W_n\); (4)
\(W_{n-1}/\operatorname{rad}W_{n-1}\) carries a positive Gram form.
Consequently \(V_{n-1}\) contains a copy of \(V_n\) with positive
quotient, and every other contribution to (T6) is controlled by (G3) or
is the identity.}

{Proof.} Parts (1) and (2) repeat the arguments of Proposition 6.1: for
\(j\le n-3\) the largest numerator index is at most \(2n-1\); for
\(j=n-2\), \(p=3\), the sole vanishing factor is \([2n]/[3]\) with
exponent one, there are no determinant zeros above the target, and the
same scalar argument kills the action on the one-dimensional radical.

For (3), the one-cap Gram matrix is
\(G_{N,1}=\delta_*I+\operatorname{Adj}(P_{N-1})\) with determinant
\([N]=[2n+2]\ne0\), and it has exactly one negative eigenvalue because
\(\delta_*\lt 2\cos(\pi/N)\) while \(\delta_*\gt2\cos(2\pi/N)\). Its
inertia is therefore \((2n,1)\) and the standard module \(W_n\) is
simple --- over the complex numbers as well, since a nonsingular real
Gram matrix stays nonsingular after complexification and the rank-one
spanning argument of Section 6.1 still applies.

\textbf{Checking the reflection hypotheses.} To avoid confusing the
chromatic variable with the quantum parameter of {[}RS{]}, write
\(\zeta=\exp(i\pi/(2n))\). The dictionary for this single application
is:

\begin{longtable}[]{@{}p{\dimexpr.42\linewidth-2\tabcolsep}p{\dimexpr.58\linewidth-2\tabcolsep}@{}}
\toprule\noalign{}
Notation in {[}RS{]} & Value here \\
\midrule\noalign{}
\endhead
\bottomrule\noalign{}
\endlastfoot
Strand number \(n\) & \(N=2n+2\) \\
Quantum parameter \(q\) & \(\zeta\) \\
Loop parameter \(\beta=q+q^{-1}\) & \(\delta_*\) \\
Smallest positive \(\ell\) with \(q^{2\ell}=1\) & \(2n\) \\
\(V_{N,1}\) and \(V_{N,2}\) & \(W_n\) and \(W_{n-1}\) \\
\end{longtable}

The critical condition of {[}RS, Section 7{]} is divisibility of
\(N-2p+1\) by \(2n\). For \(p=1,2\) these quantities are \(2n+1\) and
\(2n-1\): both labels are noncritical, they straddle the single critical
line \(N-2p+1=2n\), and \(0\lt 2-1\lt 2n\). Equivalently the
through-line labels \(2n\) and \(2n-2\) are reflected across \(2n-1\).
Since \(\delta_*\ne0,\pm2\), the exceptional parameters do not arise,
and Theorem 7.2 of {[}RS{]} gives
\(\operatorname{rad}V_{N,2}\cong L_{N,1}\) in the direction needed;
simplicity of \(V_{N,1}=W_n\) gives \(L_{N,1}=W_n\). Hence

\[ R:=\operatorname{rad}W_{n-1}\;\cong\;W_n, \qquad \dim R=N-1. \tag{G5} \]

\textbf{Counting signs before specialization.} Positivity of the
quotient remains to be proved; it cannot be inferred from the
isomorphism alone. The largest real zero of \([k]_\delta\) is
\(2\cos(\pi/k)\), and it is simple. Hence (G4) has no zero for
\(\delta\gt\delta_o:=2\cos(\pi/(2n+1))\), and the Gram matrix is
positive there by continuation from large \(\delta\). At \(\delta_o\)
its determinant has a simple zero, so precisely one eigenvalue vanishes
and changes sign; in the interval \(\delta_*\lt\delta\lt\delta_o\) there
are no further zeros, so the Gram matrix has exactly one negative
direction throughout. At \(\delta_*\) itself, \([2n+1]=-1\) and
\([2]=\delta_*\ne0\), so the determinant zero at the target has order
exactly \(N-1\).

\textbf{The derivative form.} Set \(x=\delta-\delta_*\) and choose a
fixed real basis consisting of \(R\) --- available as the kernel of the
real matrix \(G(\delta_*)\) --- and a vector-space complement. Then

\[ G(\delta_*+x)=\begin{pmatrix} x\Gamma+O(x^2) & x\Xi+O(x^2)\\[2pt] x\Xi^{T}+O(x^2) & \Theta+O(x) \end{pmatrix}, \tag{G6} \]

with \(\Theta\) nonsingular: a complement vector orthogonal to the
complement is already orthogonal to \(R\), hence lies in the radical,
hence in \(R\), hence vanishes. Writing \(A_x,\Xi_x,\Theta_x\) for the
three blocks, the Schur complement satisfies

\[ A_x-\Xi_x\Theta_x^{-1}\Xi_x^{T}=x\Gamma+O(x^2), \qquad \det G(\delta_*+x)=x^{\,N-1}\,\det\Theta_x\,\det\bigl(\Gamma+O(x)\bigr). \tag{G6a} \]

Because the order of the determinant zero is exactly \(N-1=\dim R\), the
matrix \(\Gamma\) must be nonsingular --- a conclusion that uses the
\emph{order} of vanishing, not merely the sign of the determinant. Two
preliminary observations, each one line. First, the radical is a
submodule: for \(v\in R\) and arbitrary \(w\), self-adjointness gives
\(\langle e_i(\delta_*)v,\,w\rangle=\langle v,\,e_i(\delta_*)w\rangle=0\),
so \(e_i(\delta_*)v\in R\) --- the invariance used silently below.
Second, the standard generator matrices are polynomial in \(\delta\) in
the fixed diagram basis. Differentiating the self-adjointness identity
\(G(\delta)e_i(\delta)=e_i(\delta)^{T}G(\delta)\) and sandwiching
between vectors of \(R\) kills both terms containing \(G(\delta_*)\),
and since \(e_i(\delta_*)\) preserves \(R\),

\[ \Gamma\,e_i^{R}=(e_i^{R})^{T}\Gamma, \tag{G6b} \]

where \(e_i^{R}\) is the restricted action at the target. No moving
kernel and no varying basis are needed.

Now transfer the one-cap Gram form to \(R\) through the isomorphism (G5)
and call it \(K\). The intertwiner can be taken real: its real or
imaginary part is a nonzero real intertwiner whose complexification is
injective by simplicity, hence an isomorphism by equal dimensions. Thus
\(K\) is a real invariant form of inertia \((2n,1)\). \medskip
\textbf{Lemma 6.4 (sign of the radical scalar).} \emph{In the setting above, \(\Gamma=\gamma K\) for a real scalar \(\gamma\), and \(\gamma>0\).}

\textbf{Proof.} The identities
\(\Gamma e_i^{R}=(e_i^{R})^{T}\Gamma\) and \(Ke_i^{R}=(e_i^{R})^{T}K\)
make \(K^{-1}\Gamma\) commute with all generators; complex simplicity
and Schur\textquotesingle s lemma force \(\Gamma=\gamma K\) with
\(\gamma\) a nonzero real scalar. If
\(\gamma\lt0\), the block \(x\Gamma+O(x^2)\) of the Schur congruence below would alone contribute \(2n\gt1\)
negative directions for small \(x\gt0\), contradicting the single negative direction counted
above; hence \(\gamma\gt0\). \(\blacksquare\)

\medskip
By Lemma 6.4, \(\Gamma=\gamma K\) with \(\gamma>0\). For small \(x\gt0\), Schur congruence
block-diagonalizes (G6) into \(x\Gamma+O(x^2)\) and \(\Theta+O(x)\): the first block carries exactly one
negative direction, leaving none for \(\Theta+O(x)\), and nonsingularity
then gives \(\Theta\gt0\). In particular the Gram form at the target is
positive semidefinite and induces a positive-definite form on the
quotient by \(R\). This proves (4).

\textbf{Lemma 6.5 (fusion exactness).} \emph{Applying \(F\) to the exact sequence \(0\to R\to W_{n-1}\to W_{n-1}/R\to0\) yields an exact sequence of \(F\,\mathrm{TL}_N F\)-modules.}

\textbf{Proof.} All
module maps commute with \(F\), and a quotient vector in the image of
\(F\) lifts and then projects by \(F\). This does not require \(F\) to
be central; the projected sequence consists of modules for the corner
algebra \(F\,\mathrm{TL}_N\,F\), which contains the physical updates. \(\blacksquare\)

\medskip
\textbf{Fusion and the spectrum.} By Lemma 6.5,

\[ 0\longrightarrow V_n\longrightarrow V_{n-1}\longrightarrow F\bigl(W_{n-1}/R\bigr)\longrightarrow 0, \tag{G7} \]

with positive last term of dimension \(n(n-1)/2\), while the embedded
carrier has dimension \(n\). The embedding intertwines each local
\(U_i\) and \(S_i\), hence every sweep word; in a basis extending the
embedded carrier, a sweep is block triangular, and its characteristic
polynomial is the product of those on the carrier and on the quotient
--- \emph{whether or not} (G7) splits. The quotient bound (G3) therefore
accounts for its entire spectrum, and the remaining sector \(V_{n+1}\)
is the identity. {\(\square\)}

{Corollary 6.3.} \emph{Under \(0\lt r_v,r_e\lt1\), at \(B_{2n}\) the
positive simple eigenvalue \(\lambda_*\gt1\) appears exactly once in
each of \(V_n\) and \(V_{n-1}\), and every other eigenvalue in the trace
has modulus at most \(1\). If instead only \(0\lt r_v\lt1\) and
\(r_e\gt0\), the positive subquotients have radius at most
\(\max(1,r_e)^{\,n}\), while the nonleading carrier eigenvalues still
have modulus below \(1\).}

{Remark.} Corollary 6.3 is a statement about \emph{appearances in the
standard trace}, not a claim that the algebra is semisimple at the root
of unity.

\hypertarget{s7}{%
\section{Two real-root lemmas}\label{s7}}

Only real sign changes are used in this paper. The Beraha--Kahane--Weiss
theory {[}BKW{]} and its extension by Brown and Otto {[}BO{]} are the
classical results on limiting zeros of recurrence families; their
complex-zero conclusions do not by themselves yield real roots, and the
gap is genuine: \(f_m(x)=x^2+1/m\) has complex zeros converging to the
real point \(0\) and no real zero for any \(m\). A real limit point of
zeros must be carefully distinguished from a limit of real zeros. The
two lemmas below establish the latter, through endpoint signs, under
exactly the hypotheses provided by Sections 5--6; no priority is claimed
for the general philosophy of dominant-term cancellation.

{Lemma 7.1 (vanishing single coefficient).} \emph{Suppose, near \(x=0\),
that \(A\) and \(\lambda\) are real analytic, \(A(0)=0\), \(A'(0)\ne0\),
\(\lambda(0)\gt0\), and}
\[ f_m(x)=A(x)\,\lambda(x)^{2m}+E_m(x), \qquad |E_m(x)|\le C\,R^{2m}, \]
\emph{where each \(E_m\) is real and continuous and the bound holds
uniformly on a fixed neighborhood on which \(\lambda(x)\ge L\gt R\).
Then \(f_m\) has real zeros converging to \(0\).}

{Proof.} On the two endpoints of any sufficiently small fixed interval
about \(0\), \(A\) takes opposite signs; for large \(m\) the leading
term dominates the remainder at both endpoints; the intermediate value
theorem supplies a zero; and taking closest zeros, or a nested sequence
of intervals, gives convergence. {\(\square\)}

{Lemma 7.2 (cancellation of two coefficients).} \emph{Suppose
\(A,B,\alpha,\beta\) are real analytic near \(0\), with}
\[ \alpha(0)=\beta(0)=L\gt0,\qquad A(0)+B(0)=0,\qquad c=(A+B)'(0)\ne0, \]
\emph{and \(E_m\) is real, continuous, and uniformly bounded by
\(CR^{2m}\) with \(R\lt L\) on a fixed neighborhood. Then}
\[ f_m(x)=A(x)\,\alpha(x)^{2m}+B(x)\,\beta(x)^{2m}+E_m(x) \] \emph{has a
real zero tending to \(0\) for every sufficiently large \(m\).}

{Proof.} Normalize the first two terms by \(L^{2m}\) and call the result
\(D_m\). Then
\[ D_m(0)=0,\qquad D_m'(0)=c+2md,\qquad d=\frac{A(0)\alpha'(0)+B(0)\beta'(0)}{L}. \tag{R1} \]
For all large \(m\), \(|c+2md|\ge b\gt0\): if \(d=0\) this is \(|c|\);
otherwise the expression grows. Choose \((R/L)^2\lt h\lt1\) and set
\(\varepsilon_m=h^{m}\). On \(|x|\le\varepsilon_m\), analyticity and
positivity near \(0\) give
\(|D_m''(x)|\le C_1m^2\exp(C_2m\varepsilon_m)\), so
Taylor\textquotesingle s theorem yields
\[ D_m(\pm\varepsilon_m)=\pm\varepsilon_m\,(c+2md)+O\!\bigl(m^2\varepsilon_m^2e^{C_2m\varepsilon_m}\bigr), \]
with error \(o(\varepsilon_m)\); and the remainder satisfies
\[ \frac{|E_m(\pm\varepsilon_m)|}{L^{2m}\varepsilon_m}\le C\Bigl(\frac{(R/L)^2}{h}\Bigr)^{m}\longrightarrow0 . \]
The endpoint signs therefore differ for all sufficiently large \(m\).
{\(\square\)}

{Remark.} Lemma 7.2 deliberately does not require
\(\alpha'(0)\ne\beta'(0)\); it covers identical branches. Simplicity
within each separate matrix block supplies analytic branches even when
the two agree at the target --- the even-index situation, where
coincidence of two branch values at one point must not be confused with
a repeated characteristic root throughout a neighborhood. For the
application to (T6), the remainder bound follows from spectral
continuity and finite dimension; the trace of a power is the sum of
eigenvalue powers with algebraic multiplicity, including for
nondiagonalizable matrices, so no Jordan-block factor in \(m\) is
needed. The width and the disk are fixed before the neighborhood is
chosen.

\hypertarget{spectral-localization-and-branch-bookkeeping}{%
\subsection{Spectral localization and branch
bookkeeping}\label{spectral-localization-and-branch-bookkeeping}}

The two lemmas above consume hypotheses --- a dominant analytic branch,
a uniformly exponentially smaller remainder --- that Sections 5--6
establish only \emph{at} the target. The passage from one point to a
neighborhood is elementary; we record it once, with its quantifiers in
order, rather than re-derive it informally at each use. Throughout, the
width \(n\) and the disk are fixed first; all neighborhoods and
constants may depend on them.

{Lemma 7.3 (neighborhood gap).} \emph{Let \(q\mapsto P(q)\) be a
\(d\times d\) matrix family with entries real analytic near \(q_*\).
Suppose at \(q_*\) the spectrum consists of an algebraically simple
eigenvalue \(\lambda_*\gt M\ge0\) together with eigenvalues of modulus
at most \(M\). Then there are a neighborhood \(U\ni q_*\), constants
\(L\gt R\gt M\), and a real analytic function \(\lambda\) on \(U\) with
\(\lambda(q_*)=\lambda_*\), such that for every \(q\in U\):
\(\lambda(q)\ge L\) is a simple eigenvalue of \(P(q)\), every other
eigenvalue of \(P(q)\) has modulus at most \(R\), and}
\[ \operatorname{tr}\bigl(P(q)^t\bigr)=\lambda(q)^t+E(q,t),\qquad |E(q,t)|\le(d-1)\,R^{\,t}\quad(t\ge1). \]

{Proof.} The characteristic polynomial has real analytic coefficients.
Since \(\lambda_*\) is a simple root, the implicit function theorem
gives a real analytic branch \(\lambda(q)\) of simple roots near
\(q_*\). Choose \(M\lt R\lt L\lt\lambda_*\) and apply
Rouché\textquotesingle s theorem on the fixed circle \(|z|=R\): at
\(q_*\) exactly \(d-1\) roots lie inside, the characteristic polynomial
varies uniformly on the compact circle, so the count persists on a
neighborhood; shrinking \(U\) also keeps \(\lambda(q)\ge L\). For the
trace formula, bring \(P(q)\) to Schur upper-triangular form: the trace
of a power of a triangular matrix is the sum of the \(t\)-th powers of
its diagonal, so \(\operatorname{tr}(P^t)\) is the sum of eigenvalue
\(t\)-th powers with algebraic multiplicity --- valid with no
diagonalizability assumption and with no polynomial-in-\(t\) correction.
Splitting off the simple branch gives the display. {\(\square\)}

{Lemma 7.4 (branch bookkeeping at even targets).} \emph{At
\(q_*=B_{2n}\), under the hypotheses of Corollary 6.3 with \(M=1\) (or
the certified separated bounds of Sections 9.3--9.4), \(\lambda_*\) is
an algebraically simple eigenvalue of \(P_n(q_*)\) and of
\(P_{n-1}(q_*)\) separately. Consequently there are real analytic
branches \(\beta(q)\) and \(\alpha(q)\) of eigenvalues of \(P_n(q)\) and
\(P_{n-1}(q)\) with \(\alpha(q_*)=\beta(q_*)=\lambda_*\), and on a
neighborhood of \(q_*\),}
\[ \sum_{j}[2j+1]_{\sqrt q}\operatorname{tr}\bigl(P_j(q)^{t}\bigr) \;=\; [2n-1]_{\sqrt q}\,\alpha(q)^{t}\;+\;[2n+1]_{\sqrt q}\,\beta(q)^{t}\;+\;E_t(q), \qquad |E_t(q)|\le C_J\,R^{\,t}, \]
\emph{with \(C_J=\sum_j(\dim V_j)\sup_{q\in J}|[2j+1]_{\sqrt q}|\) on a
sufficiently small compact interval \(J\) about \(q_*\),
\(R\lt\lambda_*\) as in Lemma 7.3, and \(E_t\) real and continuous. The
coefficients of \(\alpha^t\) and \(\beta^t\) are exactly \(1\) within
their sectors: no hidden multiplicity and no Jordan factor in \(t\)
occurs for the leading branches.}

{Proof.} Simplicity in \(P_n(q_*)\) is Lemma 5.1. For \(P_{n-1}(q_*)\):
by the block-triangular factorization established in the proof of
Proposition 6.2 (``Fusion and the spectrum''), its characteristic
polynomial is the product of the carrier copy\textquotesingle s and the
quotient\textquotesingle s; \(\lambda_*\) is simple in the carrier
factor by Lemma 5.1, and it is not a root of the quotient factor, whose
spectral radius is strictly below \(\lambda_*\) by (G3) under the filter
hypotheses, by the certified companion bound at \(B_8\), or by the bound
of Section 9.4 at \(B_{10}\). For the trace formula: upper-triangularize
any square matrix \(A\) over \(\mathbb C\) (Schur); then
\(\operatorname{tr}(A^{t})\) is the sum of the \(t\)-th powers of the
eigenvalues with algebraic multiplicity --- for every \(t\ge1\),
diagonalizable or not. Applying this in each sector and splitting off
the two simple leading eigenvalues, whose branches exist and stay simple
on a neighborhood by Lemma 7.3, leaves \(E_t(q)\) as a weighted sum of
eigenvalue powers. Each nonleading eigenvalue has modulus at most \(R\),
and the quantum weights are bounded on \(J\), giving the displayed
\(C_JR^t\) estimate; it is real because the total and the two leading
terms are real. For odd indices, the corresponding sum has the single
leading term with coefficient \([2n+1]_{\sqrt q}\). Multiplication of
the complete sum by \(q\) gives the coefficient \(q[2n+1]_{\sqrt q}\)
for the full normalized polynomial and multiplies the remainder bound by
\(\sup_{q\in J}q\). {\(\square\)}

{Remark.} Lemma 7.4 is what licenses taking \(A(q)=q[2n-1]_{\sqrt q}\)
and \(B(q)=q[2n+1]_{\sqrt q}\) \emph{exactly} in Lemma 7.2 --- the
quantum-integer weights multiply eigenvalue powers with multiplicity
one, not merely to leading order. Odd quantum integers are polynomials
in \(q\), so \(A,B\) are real analytic. At \(B_8\), the logical chain is
assembled in one place for the referee: the divisibility of
characteristic polynomials is the algebraic statement (G7) of
Proposition 6.2 applied at the target; the certificate bounds the
quotient\textquotesingle s companion matrix at the target; and Lemma 7.3
extends all bounds to a neighborhood.

\textbf{Corollary 7.5 (branches at odd targets).} \emph{In the setting of Proposition 6.1, at \(q_*=B_{2n+1}\) with \(0<r_v,r_e<1\): the distinguished sector carries \(\lambda_*>1\) as an algebraically simple eigenvalue, every other trace eigenvalue of every sector has modulus at most \(M<\lambda_*\) at the target, and Lemma 7.3, applied to the distinguished sector alone, yields a real analytic branch \(\lambda(q)\) with \(\lambda(q_*)=\lambda_*\), appearing in its trace with coefficient exactly \(1\); and since the sectors are finitely many and their spectra vary continuously, a common compact interval \(J\ni q_*\) can be chosen on which \(\lambda(q)\ge L\) while every other eigenvalue of every sector has modulus at most \(R\), for constants \(L>R>M\). On \(J\),}
\[ \sum_{j=0}^{n+1}[2j+1]_{\sqrt q}\operatorname{tr}\bigl(P_j(q)^t\bigr)=A_{\mathrm{odd}}(q)\,\lambda(q)^t+E(q,t),\qquad A_{\mathrm{odd}}=[2n+1]_{\sqrt q},\qquad |E(q,t)|\le C\,R^{\,t}, \]
\emph{with \(C\) bounded by the sum over all sectors of dimension times \(\sup_J|[2j+1]_{\sqrt q}|\). The amplitude \(A_{\mathrm{odd}}\) has a simple zero at \(q_*\), and the hypotheses of Lemma 7.1 are met with it; multiplication by the positive outer factor \(q\), as in (T6), changes the amplitude to \(qA_{\mathrm{odd}}\) and affects neither the zero nor the sign-change argument.}

\textbf{Proof.} Lemma 7.3 applies to the distinguished sector and yields the branch; for the remaining sectors, their eigenvalues at \(q_*\) lie in \(|z|\le M\), the eigenvalues of a finite collection of analytic matrix families vary continuously, and shrinking to a common compact \(J\) keeps all of them inside \(|z|\le R\) --- this, not Lemma 7.3 alone, is what bounds the competitors. The simple branch then appears once with its quantum-integer weight, identified exactly as \([2n+1]_{\sqrt q}\) from the bracketed sum in (T6), and every other weight--eigenvalue product is absorbed into \(E\), the weights being bounded on the compact \(J\). \(\blacksquare\)

\hypertarget{s8}{%
\section{All indices at least eleven}\label{s8}}

\hypertarget{odd-indices-1}{%
\subsection{Odd indices}\label{odd-indices-1}}

Fix \(\ell=2n+1\ge11\). Since \(B_\ell\gt B_{10}\), Proposition 2.2
supplies a fixed filter disk satisfying (D8); Lemma 5.1 and Proposition
6.1 then give a single positive simple eigenvalue \(\lambda_*\gt1\) in
the carrier \(V_n\), with every other trace eigenvalue of modulus at
most \(1\). Its coefficient inside the bracketed sum of (T6) is
\(A(q)=[2n+1]_{\sqrt q}\), which vanishes at \(q_*=B_\ell\); with
\(\theta_*=\pi/\ell\), differentiating through \(q=2+2\cos(2\theta)\)
gives

\[ A'(q_*)=\frac{\ell}{4\sin\theta_*\,\sin(2\theta_*)}\;\gt\;0. \tag{R2} \]

Taking \(t=2m\), Lemma 7.3 and Corollary 7.5 convert the
target-point spectral picture into the neighborhood hypotheses of Lemma
7.1 for the identity (T6), and real roots converge to \(B_\ell\). This
proves every odd \(\ell\ge11\).

\hypertarget{even-indices}{%
\subsection{Even indices}\label{even-indices}}

Fix \(\ell=2n\ge12\) and choose the disk the same way. Corollary 6.3
provides one simple leading eigenvalue in each of \(V_{n-1}\) and
\(V_n\), equal at the target; their coefficients are
\(A(q)=q[2n-1]_{\sqrt q}\) and \(B(q)=q[2n+1]_{\sqrt q}\), with
\(A(q_*)=q_*\) and \(B(q_*)=-q_*\). The pair cancels, and the
cancellation is transverse:

\[ A(q)+B(q)=q\sqrt q\,[2n]_{\sqrt q}, \qquad (A+B)'(q_*)=\frac{2n\,q_*\,\delta_*}{4\sin\theta_*\,\sin(2\theta_*)}\;\gt\;0. \tag{R3} \]

Lemma 7.4 supplies the analytic branches \(\alpha,\beta\), their unit
coefficients, the exact identification of \(A,B\), and the uniform
remainder bound on a neighborhood; Lemma 7.2 then applies to (T6)
divided by its positive prefactor. This proves every even \(\ell\ge12\).

{Remark.} No gap uniform in \(n\) is asserted or needed anywhere. Each
target requires one finite disk, one width, and one positive gap
depending on them --- and nothing more.

\hypertarget{s9}{%
\section{The finite set of smaller indices}\label{s9}}

\hypertarget{integer-values}{%
\subsection{Integer values}\label{integer-values}}

For \(\ell=2,3,4,6\) the targets are \(0,1,2,3\), and the single graph
\(K_4\) --- a simple sphere triangulation --- already has
\(\chi_{K_4}(q)=q(q-1)(q-2)(q-3)\); consistently with Theorem 1.1, no
accumulation claim is made at these values.

\hypertarget{the-bipyramid-at-b_5}{%
\subsection{\texorpdfstring{The bipyramid at
\(B_5\)}{The bipyramid at B\_5}}\label{the-bipyramid-at-b_5}}

Let \(D_t\) be the bipyramid over the cycle \(C_t\): two nonadjacent
poles, each adjacent to every cycle vertex --- a simple sphere
triangulation for \(t\ge3\). Splitting the count according to whether
the pole colors agree gives

\[ \chi_{D_t}(q)=q\Bigl[(q-2)^t+(q-1)(q-3)^t+(-1)^t\,(q^2-3q+1)\Bigr]. \tag{S1} \]

At \(t=2m\), in a small closed neighborhood of \(B_5\in(2,3)\), the
first two powers tend to zero uniformly, and the factor \(q^2-3q+1\) has
a simple zero at \(B_5\) with positive derivative: opposite endpoint
signs yield real roots tending to \(B_5\). At \(q=B_5\) the last term
vanishes --- \(B_5\) is a root of \(q^2-3q+1\) --- and the bracket
equals \((B_5-2)^{2m}+(B_5-1)(B_5-3)^{2m}\gt0\), a sum of even powers
with positive coefficients; so the produced roots differ from their
limit.

\hypertarget{b_7-b_8-b_9-finite-exact-certificates}{%
\subsection{\texorpdfstring{\(B_7\), \(B_8\), \(B_9\): finite exact
certificates}{B\_7, B\_8, B\_9: finite exact certificates}}\label{b_7-b_8-b_9-finite-exact-certificates}}

Below \(B_{10}\) the filter of Proposition 2.2 is unavailable
(\(\kappa\le1\)), and we use the elementary disk instead:

\[ (u,v,e)=(-1,\;q-3,\;q-2), \qquad S_i=I-U_i,\qquad \eta=1, \tag{S2} \]

with \(n=3\) at \(B_7\) and \(n=4\) at \(B_8,B_9\); Section 3 shows the
closures are simple triangulations. Here \(r_e\gt1\), so the soft
contraction argument of Proposition 6.1 is not available, and the sector
bounds are instead \emph{certified}: Appendix A constructs every
required small matrix literally from (T3) and proves, in exact
outward-rounded integer arithmetic, that a suitable power of \(10P/11\)
has sup-norm strictly below the stated rational bound, hence
\(\rho(P)\lt 11/10\):

\begin{longtable}[]{@{}lllll@{}}
\toprule\noalign{}
Target & Sector or quotient & Dimension & \(k\) & Upper bound on
\(\|(10P/11)^k\|_\infty\) \\
\midrule\noalign{}
\endhead
\bottomrule\noalign{}
\endlastfoot
\(B_7\) & \(j=0\) & 3 & 2 & \(394/1000\) \\
\(B_7\) & \(j=1\) & 6 & 4 & \(179/1000\) \\
\(B_7\) & \(j=2\) & 6 & 32 & \(264/1000\) \\
\(B_8\) & \(j=0\) & 6 & 2 & \(845/1000\) \\
\(B_8\) & \(j=1\) & 15 & 8 & \(591/1000\) \\
\(B_8\) & \(j=2\) & 15 & 64 & \(49/1000\) \\
\(B_8\) & \(V_3/V_4\), via companion matrix & 6 & 16 & \(152/1000\) \\
\(B_9\) & \(j=0\) & 6 & 8 & \(210/1000\) \\
\(B_9\) & \(j=1\) & 15 & 16 & \(347/1000\) \\
\(B_9\) & \(j=2\) & 15 & 16 & \(195/1000\) \\
\(B_9\) & \(j=3\) & 10 & 64 & \(57/1000\) \\
\end{longtable}

The top sector is the identity and needs no certificate. The carrier has
\(0\lt r_v\lt1\), so Lemma 5.1 controls all its nonleading eigenvalues.
At \(B_8\), the exact sequence (G7) identifies a second carrier copy
inside \(V_3\): its characteristic polynomial divides that of \(V_3\)
exactly, \emph{by this algebraic theorem} --- inclusion of zero in an
interval remainder would not alone establish divisibility --- and
Appendix A then encloses the quotient coefficients and bounds the
spectral radius of their companion matrix, the fourth \(B_8\) row above.
For the leading eigenvalue, the same exact arithmetic verifies

\[ \Bigl(\tfrac{13}{10}+\tau\Bigr)^{2}-4c^{2}\,\tfrac{13}{10}\;\lt\;0, \qquad c^{2}=\begin{cases}1/2,& n=3,\\[2pt] B_5/4,& n=4,\end{cases} \]

so every needed competitor has modulus below \(11/10\) while
\(\lambda_*\gt13/10\). Lemmas 7.3--7.4, with \(M=11/10\) in place of
\(1\), then convert these certified target-point bounds into
neighborhood hypotheses; Lemma 7.1 applies at \(B_7\) and \(B_9\), and
Lemma 7.2 to the paired coefficients \(q[7]\), \(q[9]\) at \(B_8\).

\hypertarget{the-golden-boundary-b_10}{%
\subsection{\texorpdfstring{The golden boundary
\(B_{10}\)}{The golden boundary B\_\{10\}}}\label{the-golden-boundary-b_10}}

Put \(\phi=(1+\sqrt5)/2\). At \(q_*=B_{10}=\phi+2\) we have \(a=\phi\)
and \(\kappa=1\): Proposition 2.2 cannot be invoked, since its proof
needs \(\kappa\gt1\) strictly. But the approximating curve (D10)
survives as \(X(s)=(1,e^s,e^{-s})\), and exact computation in
\(\mathbb Q(\phi)\) gives

\[ \mathcal R(1,2,\tfrac12)=\Bigl(-\tfrac\phi2,\;\tfrac1{2\phi^2},\;\tfrac{\phi^4}2\Bigr), \qquad \mathcal R\bigl(X(\log2)\bigr)\odot X(3\log\phi)=\tfrac\phi2\,(-1,1,1). \tag{S3} \]

This explains the choice of disk: finite fan composites approach a
boundary triple with \(r_v=r_e=1\). At that limiting triple and width
\(5\), with \(\eta_0=2/\phi\), the recurrence (H3) factors as

\[ p_5(\lambda)=(\lambda+1)\bigl((\lambda+1)^2-3\eta_0^2\lambda\bigr)\bigl((\lambda+1)^2-\eta_0^2\lambda\bigr), \]

and since \(3\eta_0^2\gt4\) while \(\eta_0^2\lt4\), a simple positive
eigenvalue is already separated from everything else \emph{at the
boundary itself}; strict local contraction of both ratios is
unnecessary. For a concrete finite disk, take

\[ R=\mathcal R(F_{21}),\qquad N_c=6480,\qquad M_c=13496,\qquad D=R^{M_c}\,\mathcal R\bigl(R^{N_c}\bigr): \tag{S4} \]

compose \(N_c\) copies, rotate the whole, compose \(M_c\) further
copies. The disk has \(21(N_c+M_c)-1=419{,}495\) internal vertices and
is simple, triangulated, and free of boundary diagonals by Section 2.
The scalar certificate of Appendix A, starting from an integer enclosure
of \(\sqrt5\), proves

\[ u\lt0,\qquad \frac{999}{1000}\lt r_v\lt 1\lt r_e\lt\frac{1001}{1000},\qquad \eta\gt\frac{123}{100}. \tag{S5} \]

Proposition 6.2 applies at \(n=5\); every positive subquotient has sweep
norm at most \((1001/1000)^5\lt101/100\); the one-dimensional radical
acts as the identity; each carrier copy has exactly one expanding
eigenvalue by Lemma 5.1; and at \(c^2=3/4\),

\[ (2+r_v)^2-6\eta^2 \;\lt\; 9-6\Bigl(\tfrac{123}{100}\Bigr)^{2}\;\lt\;0, \]

so \(\lambda_*\gt2\). There are exactly two appearances of the leading
eigenvalue, all other trace eigenvalues have modulus below \(101/100\),
and the coefficients \(q[9]\), \(q[11]\) sum to \(q\sqrt q\,[10]\) with
a simple zero at the target: Lemmas 7.3--7.4 (with \(M=101/100\)) and
then Lemma 7.2 apply. The resulting graphs have
\(|V(G_{2m})|=2+4{,}194{,}960\,m\) vertices; no claim of efficiency or
minimality is made, and the construction is specified by a finite
grammar --- listing millions of vertices is not part of the certificate.

\hypertarget{s10}{%
\section{Distinct roots and completion of the proof}\label{s10}}

It remains to show the roots produced above differ from their targets.
First, for any connected loopless graph, \(\chi_G\) is nonzero on
\((0,1)\) with sign \((-1)^{|V|-1}\): by deletion--contraction, a bridge
contributes \(\chi_G=(q-1)\chi_{G/e}\), while for a nonbridge the
deletion stays connected and the deletion and negated contraction terms
share the asserted sign; parallel edges arising from contraction are
harmless; the one-vertex case starts the induction.

For odd \(\ell=2n+1\ge5\), write \(B_\ell=2+\zeta+\zeta^{-1}\) with
\(\zeta=e^{2\pi i/\ell}\). Since \(\gcd(n,\ell)=1\), the algebraic
conjugate \(2+\zeta^{n}+\zeta^{-n}=2-2\cos(\pi/\ell)\) lies in
\((0,1)\); a chromatic root is an algebraic integer whose conjugates are
all chromatic roots of the same polynomial, so \(B_\ell\) cannot be one.
For even \(\ell=2n\) with \(n\ne5\), take \(k=n-1\) if \(n\) is even and
\(k=n-2\) if \(n\) is odd; for \(n=4\) or \(n\ge6\) one checks
\(\gcd(k,2n)=1\) and \(1/3\lt k/(2n)\lt1/2\), so the conjugate
\(2+2\cos(2\pi k/(2n))\) again lies in \((0,1)\).

The one stubborn case is \(B_{10}\), where the conjugate argument fails
--- and where Tutte\textquotesingle s golden identity, in the form
stated and derived by Fendley and Krushkal {[}FK2{]}, does exactly what
is needed: for a plane triangulation \(T\),

\[ \chi_T(\phi+2)=(\phi+2)\,\phi^{\,3|V(T)|-10}\,\chi_T(\phi+1)^{2}\;\gt\;0, \tag{S7} \]

strictly, because \(\phi+1=B_5\) is not a chromatic root by the
conjugate argument just given. So every noninteger target is approached
by roots different from itself, and a sequence of such roots contains
infinitely many distinct values.

Sections 8 and 9 together cover every \(\ell\ge11\) and every
\(2\le\ell\le10\); this completes the proof of Theorem 1.1. If one adds
the convention \(B_1=4\), choosing for each large \(\ell\) a root within
\((4-B_\ell)/2\) of \(B_\ell\) gives roots below \(4\) tending to \(4\);
as noted in Section 1, that endpoint is already Royle\textquotesingle s
theorem {[}R{]}.

\hypertarget{s11}{%
\section{Exact certificates and their logical role}\label{s11}}

It is worth being precise about what is computed and what is proved. The
finite computation is confined to the explicit scalar disk at \(B_{10}\)
and the small matrices at \(B_7,B_8,B_9\); the all-width Gram and
carrier arguments of Sections 5--6 are proofs, not extrapolations from
samples.

The logical status of these certificates is that of a computer-assisted
proof in the tradition that began, for this very
problem\textquotesingle s ancestor, with the Four Color Theorem:
finitely many inequalities about finitely many explicitly displayed
integers, verified in exact outward-rounded arithmetic, with complete
and self-contained source supplied as the supplementary files \texttt{beraha\_certificates.py}, \texttt{appendixB.py}, \texttt{beraha\_trace\_audit.py}, \texttt{beraha\_symbolic\_audit.py}, \texttt{crosscheck\_t6.py} and \texttt{crosscheck\_composite.py} (with SHA-256 manifest). The reproduction package (the six scripts, their run logs and a SHA-256 manifest) is archived at Zenodo, \href{https://doi.org/10.5281/zenodo.23371280}{doi:10.5281/zenodo.23371280}, and is also included with the arXiv submission as ancillary files. The primary finite certificates and the independent convention and cross-check audits are supplied as executable supplementary programs, with recorded run logs and SHA-256 hashes, so that verification is running short programs rather than trusting the author. Nothing asymptotic, numerical, or floating-point is
certified anywhere.

Appendix A uses a fixed-point scale \(S=2^{320}\). A scalar interval
\([l,u]\) represents \([l/S,u/S]\); products take the minimum and
maximum of the four integer endpoint products, rounded outward after
division by \(S\); division encloses a reciprocal only when the
denominator interval excludes zero. A matrix ball \((M,r)\) encloses all
matrices whose entries differ from \(M/S\) by at most \(r/S\); for the
product of two balls, the entrywise error before rescaling is bounded by
\(r_B\|M_A\|_\infty+r_A\|M_B\|_1+k\,r_Ar_B\) with \(k\) the inner
dimension, integer multiplication produces the exact midpoint numerator,
and outward rounding adds at most one unit. The row-sum bound
\((\|M\|_\infty+dr)/S\) then converts a certified power of \(10P/11\)
below \(1\) into the strict inequality \(\rho(P)\lt11/10\). The
enclosures of \(B_7,B_8,B_9\) use Machin\textquotesingle s identity,
alternating rational arctangent series, and alternating cosine bounds
within their monotone ranges; \(\sqrt q\) uses integer square roots; the
paired-sector characteristic polynomial is enclosed by
Faddeev--LeVerrier, with exact divisibility supplied by (G7) before
interval division of polynomials. Every accepted bound is a comparison
of integers; decimal printouts are informational. The recorded run
completes \(112\) integer matrix products and passes every stated bound,
writing full endpoints to JSON for inspection. No convergence or error
analysis enters: the disk is a single finite, fixed object, and every
inequality asserted about it is checked exactly.

Appendix B is a separate audit of conventions. It computes
\(\chi_G(q)=\sum_k a_k (q)_k\) by counting partitions of graph vertices
into independent sets, and compares the complete symbolic polynomial
against the fused trace in \(\mathbb Z[\delta,\delta^{-1}]\), for nine
closures of \(8\)--\(14\) vertices including nonuniform cells built from
rotated and serially composed fans; it also checks face incidences,
vertex links, connectedness and the Euler characteristic. These finite
identities support the normalization and gluing conventions; the general
justification remains the diagrammatic proof of Section 4.

For the smallest widths, the uniform fan closures are close cousins of
the recursive families with which Beraha and Kahane first exhibited
complex chromatic roots near \(4\); what Sections 5--7 add, and what
those families lack, is a dominance structure valid for all widths at
once, together with realness.

\hypertarget{s12}{%
\section{Discussion and relation to prior work}\label{s12}}

One way to summarize the mechanism: the single angle \(\theta=\pi/\ell\)
controls both the quantum-integer coefficients in the trace and the
locations of Gram degeneracies. At odd indices a coefficient vanishes;
at even indices a reflected radical plants a second copy of the carrier
whose coefficient cancels the first. That is the precise role the
trigonometric parameter plays here --- no more, and no less.

The contributions of this paper are: the disk realization of Proposition
2.2; the all-width use of Gram inertia, including the derivative
argument (G6a)--(G6b); the carrier dominance lemma; and their assembly
through the real cancellation lemmas into the full index range. The
representation-theoretic determinant, the reflection theorem, the
chromatic map, the quantum trace, and the golden identity are credited
inputs, not novelty claims; nor is the general contrast between
vanishing amplitudes and exact inter-sector coincidences claimed as new
--- see {[}BCP{]}.

\hypertarget{literature-comparison}{%
\subsection{Literature comparison}\label{literature-comparison}}

\begin{longtable}[]{@{}p{\dimexpr.23\linewidth-2\tabcolsep}p{\dimexpr.38\linewidth-2\tabcolsep}p{\dimexpr.39\linewidth-2\tabcolsep}@{}}
\toprule\noalign{}
Prior work & Established framework or result & Supplied by the present
paper \\
\midrule\noalign{}
\endhead
\bottomrule\noalign{}
\endlastfoot
Saleur {[}S{]} (1990); Kauffman--Saleur {[}KS{]} (1993) & Quantum-group
and transfer-matrix approaches, with explicitly acknowledged heuristic
steps in the accumulation argument. & Dominance for an explicit finite
disk construction (Proposition 2.2, Section 6). \\
Jacobsen--Saleur {[}BCP{]} (2008) & Boundary chromatic polynomials;
amplitude zeros and exact sector coincidences at finite width, explained
by representation theory. & Dominance in the complete closure and real
zeros within the triangulation class; exact coincidence alone is
insufficient. \\
Biggs {[}B1{]} (2001); {[}B2{]} (2004) & Compatibility-matrix traces;
decompositions for bracelet families via Specht modules. & The planar
construction, its complete closure, and bounds on all contributing
sectors. \\
Beraha--Kahane--Weiss {[}BKW{]}; Brown--Otto {[}BO{]} & Limiting-zero
criteria for recurrence families; {[}BO, Thm 2.1{]} permits
index-polynomial coefficients. & Real sign changes under the verified
hypotheses (Lemmas 7.1--7.4); branch coincidence at one point is not a
repeated root in a neighborhood. \\
Fendley--Krushkal {[}FK1, FK2{]}; Jones {[}J{]} & Chromatic-to-TL
correspondence, chromatic identities, standard quantum trace. & The
explicit derivation (T2a)--(T6) with its normalization and boundary
conventions. \\
Ridout--Saint-Aubin {[}RS{]} & Standard-module Gram determinant and
reflection structure. & The radical, fusion and derivative arguments
yielding the sector bounds of Section 6. \\
Agol--Krushkal {[}AK{]} & Upper bound on chromatic values of planar
triangulations at Beraha numbers via TL positivity (indexing shifted by
one relative to ours). & Existence and realness of nearby zeros, and
dominance; a value estimate does not imply a nearby real zero. \\
Royle {[}R{]}; Perrett--Thomassen {[}PT{]} & Triangulation roots
approaching \(4\); density for general planar graphs. & Localization at
each prescribed finite target within the triangulation class. \\
Jacobsen--Salas--Sokal {[}JSS{]} & Triangular-strip transfer-matrix
analysis of the chromatic polynomial. & The pole-pinned path geometry
here differs from strip boundary conditions; no result of {[}JSS{]} is
imported. \\
\end{longtable}

In {[}S, Section 11{]} the author himself lists the unresolved rigor ---
Bethe-state completeness, conformal spectral asymptotics, regularity in
parameter space --- and the note added in proof advances the
chromatic-line argument and mentions free boundaries; it has been
included in the comparison and does not provide the finite-disk and
Gram-inertia argument given here.

\hypertarget{scope-of-the-proposed-mechanism}{%
\subsection{Scope of the proposed
mechanism}\label{scope-of-the-proposed-mechanism}}

The filter construction changes the local disk instead of keeping a
fixed lattice cell. Its local bounds must then be combined with the
complete trace, carrier dominance, and Gram inertia. In particular, a
local bound larger than one for an elementary cell does not by itself
prove that a competing global sector dominates. The fused construction
divides only by the nonzero loop parameter and uses standard modules at
the target; it does not require singular higher projectors. These are
the specific mathematical claims to be checked against prior work and
under independent review.

\textbf{Use of AI tools.}
In preparing this work the author used large language models (Anthropic's Claude and
OpenAI's ChatGPT) as assistants for the following tasks: writing and debugging the
verification and audit scripts reproduced in Appendices~A--C; editing the exposition to
standard mathematical English and converting the manuscript to \LaTeX; searching for and
checking references; and line-by-line checking of the proofs for gaps and errors. The
mathematical ideas, the constructions, and the proofs are the author's; every argument and
every script was read, verified, and is vouched for by the author, who takes full
responsibility for the entire content of the paper, including any errors. The two
implementations compared in Appendix~C were written in separate sessions with no shared
code; ``independent'' there refers to the code, not to an independent human verifier. No
AI output was accepted as evidence for any mathematical claim: the certificates are exact
integer computations whose logs are included, and the proofs stand on the text.

\hypertarget{refs}{%
\section*{References}\label{refs}}
\addcontentsline{toc}{section}{References}

{[}Be{]} S. Beraha, \emph{Infinite non-trivial families of maps and
chromials}, Ph.D. thesis, Johns Hopkins University, 1975.

{[}JS{]} J. L. Jacobsen and J. Salas, \emph{A generalized Beraha
conjecture for non-planar graphs}, Nucl. Phys. B 875 (2013), 678--718.
arXiv:1303.5210.

{[}FK1{]} P. Fendley and V. Krushkal, \emph{Link invariants, the
chromatic polynomial and the Potts model}, Adv. Theor. Math. Phys. 14
(2010), no. 2, 507--540; esp. Definition 6.1, Lemmas 6.2, 6.4, 7.1.
arXiv:0806.3484.

{[}FK2{]} P. Fendley and V. Krushkal, \emph{Tutte chromatic identities
from the Temperley--Lieb algebra}, Geom. Topol. 13 (2009), 709--741.

{[}RS{]} D. Ridout and Y. Saint-Aubin, \emph{Standard modules, induction
and the structure of the Temperley--Lieb algebra}, Adv. Theor. Math.
Phys. 18 (2014), no. 5, 957--1041; esp. Theorems 4.7 and 7.2.
arXiv:1204.4505.

{[}J{]} V. F. R. Jones, \emph{Index for subfactors}, Invent. Math. 72
(1983), 1--25: the Markov trace on the Temperley--Lieb algebra and its
weights.

{[}S{]} H. Saleur, \emph{Zeroes of chromatic polynomials: a new approach
to the Beraha conjecture using quantum groups}, Comm. Math. Phys. 132
(1990), 657--679; Sections 8, 11 and the note added in proof.

{[}KS{]} L. H. Kauffman and H. Saleur, \emph{An algebraic approach to
the planar coloring problem}, Comm. Math. Phys. 152 (1993), 565--590,
esp. Section VIII.

{[}BCP{]} J. L. Jacobsen and H. Saleur, \emph{Boundary chromatic
polynomial}, J. Stat. Phys. 132 (2008), 707--719, esp. Sections 4--5.
arXiv:0803.2665.

{[}B1{]} N. Biggs, \emph{A matrix method for chromatic polynomials}, J.
Combin. Theory Ser. B 82 (2001), 19--29.

{[}B2{]} N. Biggs, \emph{Specht modules and chromatic polynomials}, J.
Combin. Theory Ser. B 92 (2004), 359--377.

{[}BKW{]} S. Beraha, J. Kahane and N. J. Weiss, \emph{Limits of zeroes
of recursively defined polynomials}, Proc. Nat. Acad. Sci. USA 72
(1975), 4209; and \emph{Limits of zeros of recursively defined families
of polynomials}, in: Studies in Foundations and Combinatorics, Adv. in
Math. Suppl. Stud. 1, Academic Press (1978), 213--232.

{[}BO{]} J. Brown and P. T. Otto, \emph{A generalization of the
Beraha--Kahane--Weiss theorem with graph polynomial applications},
Contrib. Discrete Math. 20 (2025), no. 2, DOI
10.55016/ojs/cdm.v20i2.76682; esp. Theorem 2.1. arXiv:2008.01784.

{[}AK{]} I. Agol and V. Krushkal, \emph{Tutte relations, TQFT, and
planarity of cubic graphs}, Illinois J. Math. 60 (2016), no. 1,
273--288; esp. Theorem 4.1. arXiv:1512.07339.

{[}R{]} G. Royle, \emph{Planar triangulations with real chromatic roots
arbitrarily close to four}, Ann. Comb. 12 (2008), 195--210,
arXiv:math/0511304; and \emph{The Beraha number \(B_{10}\) is a
chromatic root}, arXiv:1611.01941.

{[}PT{]} T. J. Perrett and C. Thomassen, \emph{Density of chromatic
roots in minor-closed graph families}, Combin. Probab. Comput. 27
(2018).

{[}P{]} T. J. Perrett, \emph{Roots of the chromatic polynomial}, Ph.D.
thesis, Technical University of Denmark, 2016, Section 6.1.

{[}JSS{]} J. L. Jacobsen, J. Salas and A. D. Sokal, \emph{Transfer
matrices and partition-function zeros for antiferromagnetic Potts
models. III. Triangular-lattice chromatic polynomial}, J. Stat. Phys.
112 (2003), 921--1017. arXiv:cond-mat/0204587.

{[}LK{]} S. J. Lomonaco, Jr. and L. H. Kauffman, \emph{Topological
quantum computing and the Jones polynomial}, preprint (2006), Theorem
5.2 (uniqueness of the normalized Markov trace).
\href{https://arxiv.org/abs/quant-ph/0605004}{arXiv:quant-ph/0605004}.

{[}M{]} R. Mactabi, \emph{Exact amplitude identities at Beraha numbers
for closed planar triangulations: certified limits at B7, B9, B10 and a
two-sided structure}, author manuscript, 19 August 2026.

\hypertarget{appA}{%
\section*{Appendix A. Executable exact finite certificates}\label{appA}}
\addcontentsline{toc}{section}{Appendix A. Executable exact finite
certificates}

The executable certificate is supplied as the supplementary file \texttt{beraha\_certificates.py}; run it with \texttt{python3}, without \texttt{-O}. Python~3 and NumPy are required, SciPy is not. The two JSON outputs are written
next to the script. The algebraic justification of the \(B_8\) quotient
is Proposition~6.2, not the interval remainder test alone.

\hypertarget{appB}{%
\section*{Appendix B. Convention audit on complete graph
polynomials}\label{appB}}
\addcontentsline{toc}{section}{Appendix B. Convention audit on complete
graph polynomials}

An independent check of the gluing and normalization conventions: the
script computes complete symbolic chromatic polynomials by counting
partitions into independent sets, builds the fused trace in
\(\mathbb Z[\delta,\delta^{-1}]\), and asserts exact equality on nine
closures, including nonuniform cells built from rotated and serially
composed fans.

\hypertarget{appC}{%
\section*{Appendix C. Reproducible supplementary audits}\label{appC}}
\addcontentsline{toc}{section}{Appendix C. Reproducible supplementary
audits}

The two cross-check scripts supplied with the author revision provide a
second graph-counting and exact-module implementation. Their development
independence has not been independently established. They use Python 3
standard-library arithmetic; unused SymPy imports and a machine-specific
path have been removed, and comparisons now assert success rather than
merely printing a failure or exiting early. The complete scripts, and the run logs of every supplementary
script, are supplied with the manuscript in the supplementary archive
(directory \texttt{logs/}).

\textbf{C.1. Uniform closures.} \texttt{crosscheck\_t6.py} passed all
nine exact integer comparisons, with \(n,t\) equal to
\((2,4),(3,4),(2,6)\) and \(q=5,6,7\).

\textbf{C.2. Composite cells.} \texttt{crosscheck\_composite.py} passed
all eight exact comparisons using rotated and serially composed fans and
the uniform \((n,t)=(4,3)\) closure, at \(q=5,7\).

\textbf{C.3. Existing certificates.} The unchanged Appendix A source
passed its stated assertions with 112 integer matrix products; the
unchanged Appendix B source passed all nine full symbolic polynomial
identities.

\textbf{C.4. Trace and local cubic audit.}
\texttt{beraha\_trace\_audit.py} checks 392 exact coloring/diagram trace
equalities on every diagram basis element for strand numbers 1 through 6
at \(q=4,9\), checks the adjacent equalizer/resummation relation, and
evaluates the three fused \(\mathrm{TL}_4\) local scalars symbolically.
These checks are finite regression tests, not proofs for all widths.

The earlier numerical root-march and Gram-spectrum tables are omitted
because their drivers and the full fixed-disk data were not supplied.
None is needed for the argument or the exact certificates above.

\end{document}